\documentclass[11pt]{article}

\usepackage{amssymb}
\usepackage{amsmath}
\usepackage{amsfonts}
\usepackage[all]{xy}
\usepackage{graphicx}
\usepackage{epsfig}
\usepackage{makeidx}
\usepackage{color}
\input{epsf.tex}
\input xy
\xyoption{all}

\newtheorem{defn}{\noindent \bf{Definition}}[section]
\newtheorem{thm} {\noindent \bf{Theorem}}[section]
\newtheorem{cor} {\noindent \bf{Corollary}}[section]

\newtheorem{remark}{\noindent \bf{Remark}}[section]
\newtheorem{lem}{\noindent \bf{Lemma}}[section]
\newtheorem{prop}{\noindent \bf{Proposition}}[section]
\numberwithin{equation}{section}
\newenvironment{prooof}[1]{\noindent$\textbf{#1. }$}
{\hspace{\fill}$_\blacksquare$}

\begin{document}

\title{\sc Equilibrium balking strategies for a clearing queueing
system in alternating environment}
\author{Antonis Economou and Athanasia Manou \\
        aeconom@math.uoa.gr and amanou@math.uoa.gr\\
        University of Athens, Department of Mathematics\\
        Panepistemiopolis, Athens 15784, Greece}
\date{\today}
\maketitle

\noindent \textbf{Abstract:} We consider a Markovian clearing
queueing system, where the customers are accumulated according to a
Poisson arrival process and the server removes all present customers
at the completion epochs of exponential service cycles. This system
may represent the visits of a transportation facility with unlimited
capacity at a certain station. The system evolves in an alternating
environment that influences the arrival and the service rates. We
assume that the arriving customers decide whether to join the system
or balk, based on a natural linear reward-cost structure. We study
the balking behavior of the customers and derive the corresponding
Nash equilibrium strategies under various levels of information.

\vspace{0.5cm}


\vspace{1cm}

\section{Introduction}

Queueing systems with batch services are often used to represent the
visits of a transportation facility at a certain station. This
allows the quantification of the congestion of the station and can
be used to take control measures (e.g. changing the frequency of the
visits), so that the quality of service is kept within acceptable
limits. The capacity of the facility is usually assumed unlimited.
This is justified, because in most applications the capacity of the
facility is chosen large enough, so that the probability that some
waiting customers cannot be accommodated is negligibly small.
Moreover, the waiting customers that cannot be served at a visit of
the facility are not in general willing to wait for its next visit
and abandon the system. Therefore, it is realistic to assume that
all present customers are removed at the visit points of the
facility. Such systems are referred to as stochastic clearing systems.\\

Stochastic clearing systems have been studied extensively in the
literature (see e.g. Stidham (1974), Serfozo and Stidham (1978),
Artalejo and Gomez-Corral (1998) and Yang, Kim and Chae (2002)).
They have been also studied in the framework of stochastic systems
subject to (total) catastrophes or disasters, where catastrophic
events are assumed to remove all the customers/units of the
system/population (see e.g. Kyriakidis (1994), Economou and Fakinos
(2003,2008), Stirzaker (2006,2007) and Gani and Swift (2007)). In
the majority of such studies the interest of the investigators lies
in the transient and/or the stationary distribution of the process
of interest. However, optimization issues for this class of systems
have also attracted the interest in the literature (see e.g.
Kyriakidis
(1999a,b), Economou (2003), Kyriakidis and Dimitrakos (2005)).\\

During the last decades, there is an emerging tendency to study
queueing systems from an economic viewpoint. In the context of
stochastic clearing systems, the optimization questions that have
been considered so far concern the central planning of the systems.
In these studies, the objective is the determination of optimal
policies for the server, about when he should remove the customers
from the system (see e.g. Stidham (1977), Kim and Seila (1993),
Economou (2003), Kyriakidis and Dimitrakos (2005)). However, to the
best of our knowledge, there are no economic studies that concern
the behavior of the customers when they are free to make decisions
to maximize their own benefit. Such considerations lead to a game
theoretic economic analysis of their behavior in
the system.\\

In general, the economic analysis of customer behavior in a queueing
system is based on some reward-cost structure which is imposed on
the system and reflects the customers' desire for service and their
unwillingness to wait. Customers are allowed to make decisions about
their actions in the system, for example they may decide whether to
join or balk, to wait or abandon, to retry or not etc. The customers
want to maximize their benefit, taking into account that the other
customers have the same objective, and so the situation can be
considered as a game among them. In this type of studies, the main
goal is to find individual and social optimal strategies. The study
of queueing systems under a game-theoretic perspective was initiated
by Naor (1969) who studied the $M/M/1$ model with a linear
reward-cost structure. Naor (1969) assumed that an arriving customer
observes the number of customers and then makes his decision whether
to join or balk (observable case). His study was complemented by
Edelson and Hildebrand (1975) who considered the same queueing
system but assumed that the customers make their decisions without
being informed about the state of the system. Since then, there is a
growing number of papers that deal with the economic analysis of the
balking behavior of customers in variants of the $M/M/1$ queue, see
e.g. Hassin and Haviv (1997) ($M/M/1$ queue with priorities),
Burnetas and Economou (2007) ($M/M/1$ queue with setup times), Guo
and Zipkin (2007) ($M/M/1$ queue with various levels of information
and non-linear reward-cost structure), Hassin (2007) ($M/M/1$ queue
with various levels of information and uncertainty in the system
parameters), Economou and Kanta (2008a,b) ($M/M/1$ queue with
compartmented waiting space, $M/M/1$ queue with unreliable server),
Sun, Guo and Tian (2010) ($M/M/1$ queue with setup/closedown times)
and Zhang and Wang (2010) ($M/M/1$ queue with delayed repairs). The
monographs of Hassin and Haviv (2003) and Stidham (2009) summarize
the main approaches and several results in the broader area of the
economic analysis of queueing systems.\\

The aim of the present paper is to study the equilibrium behavior of
the customers regarding balking in the framework of a Markovian
clearing queueing model. The balking behavior of customers in
stochastic clearing systems that model transportation stations is
important and should be taken into account if one wants to obtain a
reliable representation of what is going on in these systems.
However, such systems usually evolve in random environment, i.e.
there is some external process that influences the arrival and the
service rates. In the present study we will concentrate on a
clearing system evolving in an alternating random environment
(modeled by a 2-state continuous-time Markov chain). We will
determine equilibrium balking strategies for the customers under
various levels of information. In particular, we will consider
several information cases, as an arriving customer may observe or
not the number of customers in the system and/or the state of the
environment, before making his
decision about whether to join or balk.\\

The paper is organized as follows. In Section 2, we describe the stochastic dynamics of the model, the reward-cost structure and the decision framework (information cases). In Section 3 we consider those cases, where the strategies of the other customers do not influence the expected net benefit of a tagged customer. These are the unobservable cases, where the tagged customer does not observe the number of customers in the system before making his decision and the fully observable case, where he observes both the number of customers in the system and the state of the environment. In all these cases, we show that the expected net benefit of a tagged customer depends only on his strategy and not on the strategies followed by the other customers, a fact that implies the existence of dominant strategies. This is a special feature of the system that is related to the nature of the stochastic clearing mechanism.
In Sections 4 and 5 we consider the almost observable case, where the customers get informed upon arrival about the number of waiting customers in the station but not about the state of the environment. In this case, the waiting customers do not imply any additional cost on the individual, but their presence provides a signal about the clearing rate. Depending on the parameters of the model, a large number of waiting customers may increase or decrease the conditional probability that the clearing rate is the slow one. In Section 4, we present some preliminary results. More concretely, we first compute the stationary distributions of the system, when the customers follow either a threshold or a reverse-threshold strategy. Then we compute the net benefit of an arriving customer who decides to join, given that he observes $n$ customers and that the others follow a threshold or a reverse-threshold strategy. In Section 5, we conclude our study and we characterize all equilibrium strategies within the class of threshold and
reverse-threshold strategies. The main contribution of the paper is an algorithm that computes efficiently all these equilibrium strategies.  In Section 6, we summarize our findings and discuss the Follow-The-Crowd and Avoid-The-Crowd notions for this model as well as the problem of social optimization.\\

\section{The model}

We consider a transportation station with infinite waiting space
that operates in an alternating environment. The environment is
specified by a 2-state continuous-time Markov chain $\{E(t)\}$, with
state space $S^E=\{1,2\}$ and transition rates $q_{ee'}$, for $e\ne
e'$. Whenever the environment is at state $e$, customers arrive
according to a Poisson process at rate $\lambda_e$, whereas a
transportation facility visits the station according to a Poisson
process at rate $\mu_e$. The two Poisson processes are assumed
independent. At the visit epochs of the transportation facility all
customers are served instantaneously and removed from the station.
Therefore, we have a stochastic clearing
system in an alternating random environment.\\

We represent the state of the station at time $t$ by a pair
$\left(N(t),E(t)\right)$, where $N(t)$ records the number of
customers at the station and $E(t)$ denotes the environmental state.
The stochastic process
$\left\{\left(N(t),E(t)\right):t\geq0\right\}$ is a continuous-time
Markov chain with state space $S^{N,E}=\{(n,e):\; n\geq 0,\;
e=1,2\}$ and its non-zero transition rates are given by
\begin{eqnarray}
q_{(n,e)(n+1,e)}&=&\lambda_e,\; n\geq 0,\; e=1,2,\label{TranRate1}\\
q_{(n,e)(0,e)}&=&\mu_e,\; n\geq 1,\; e=1,2,\label{TranRate2}\\
q_{(n,1)(n,2)}&=&q_{12},\; n\geq 0,\label{TranRate3}\\
q_{(n,2)(n,1)}&=&q_{21},\; n\geq 0.\label{TranRate4}
\end{eqnarray}
We define $\rho_e=\frac{\lambda_e}{\mu_e}$, $e=1,2$. The value of $\rho_e$ can be thought of as a measure of congestion of the system under the environmental state $e$, as it expresses the mean number of customers accumulated between two successive visits of the transportation facility (given that the environment remains continuously in state $e$).

We are interested in the behavior of customers, when they have the
option to decide whether to join or balk. We assume that a customer
receives a reward of $R$ utility units for completing service.
Moreover, a customer accumulates costs at a rate of $C$ utility
units per time unit that he remains in the system. We also assume
that customers are risk neutral and wish to maximize their net
benefit. Finally, their decisions are assumed irrevocable, in the
sense that neither reneging of entering
customers nor retrials of balking customers are allowed.\\

Since all customers are assumed indistinguishable, we can consider
the situation as a symmetric game among them. Denote the common set
of strategies (set of available actions) and the utility (payoff) function by
$\mathcal{S}$ and $\mathcal{U}$ respectively. More concretely, let $\mathcal{U}\left(s_{tagged},
s_{others}\right)$ be the payoff of a tagged customer who follows
strategy $s_{tagged}$, when all other customers follow $s_{others}$.
A strategy $s_{1}$ is said to dominate strategy $s_{2}$ if
$\mathcal{U}\left(s_{1}, s\right)\geq \mathcal{U}\left(s_{2}, s\right)$, for every $s
\in \mathcal{S}$. A strategy $s_{*}$ is said to be dominant if it
dominates all other strategies in $\mathcal{S}$. A strategy $\tilde{s}$ is
said to be a best response against a strategy $s_{others}$, if
$\mathcal{U}\left(\tilde{s}, s_{others}\right)\geq
\mathcal{U}\left(s_{tagged},s_{others}\right)$, for every $s_{tagged}\in \mathcal{S}$.
Finally, a strategy $s_{e}$ is said to be a (symmetric) Nash
equilibrium, if and only if it is a best response against itself,
i.e. $\mathcal{U}\left(s_{e}, s_{e}\right)\geq \mathcal{U}\left(s, s_{e}\right)$, for
every $s\in \mathcal{S}$. The intuitive interpretation of a Nash equilibrium
is that it is a stable point of the game, in the sense that if all
customers agree to follow it, then no one can benefit by deviating from
it. We remark that the notion of a dominant strategy is stronger
than the notion of an equilibrium. In fact, every dominant strategy
is an equilibrium, but the converse is not true. Moreover, while
equilibrium strategies exist in most situations,
dominant strategies rarely do.\\

In the next sections we obtain customer equilibrium strategies for
joining/balking. We distinguish four cases depending on the
information available to the customers at their arrival instants,
before the decision is made:
\begin{itemize}
\item Fully unobservable case: Customers do not observe $N(t)$
nor $E(t)$.
\item Almost unobservable case: Customers do not observe $N(t)$,
but observe $E(t)$.
\item Fully observable case: Customers observe both $N(t)$ and
$E(t)$.
\item Almost observable case: Customers observe $N(t)$, but do not observe $E(t)$.
\end{itemize}

From a methodological point of view, the first three cases are similar and they lead to dominant strategies, so we study all of them in Section \ref{Unobservable-Section}. The almost observable case which is the most interesting and methodologically demanding is treated in Sections \ref{Observable-Section-Preliminaries} and
\ref{Observable-Section-Equilibrium}.

\section{The unobservable and the fully observable cases: Dominant strategies}\label{Unobservable-Section}

Let $S_{e}$ denote the time till the next arrival of the transportation facility, given that the environment is at state $e$.
A moment of reflection shows that $S_{e}$ is independent of the number of customers in the system, because of the mechanism of the total
removals of customers at the visits of the facility and the
memoryless property of the exponential distribution. By employing a
first-step argument, conditioning on the next transition of the
Markov chain $\{(N(t),E(t))\}$ that is either a visit of the
facility or a change in the environment, we obtain the equations

\begin{eqnarray}
E[S_{1}]&=&\frac{1}{\mu_1+q_{12}}+\frac{\mu_1}{\mu_1+q_{12}}\cdot 0+\frac{q_{12}}{\mu_1+q_{12}} \cdot E[S_{2}],\label{ES1-equation}\\
E[S_{2}]&=&\frac{1}{\mu_2+q_{21}}+\frac{\mu_2}{\mu_2+q_{21}}\cdot
0+\frac{q_{21}}{\mu_2+q_{21}} \cdot E[S_{1}].\label{ES2-equation}
\end{eqnarray}

Solving \eqref{ES1-equation}-\eqref{ES2-equation} yields
\begin{eqnarray}
E[S_{1}]&=&\frac{\mu_2+q_{21}+q_{12}}{\mu_1 \mu_2+\mu_1 q_{21}+\mu_2 q_{12}},\label{ES1}\\
E[S_{2}]&=&\frac{\mu_1+q_{21}+q_{12}}{\mu_1 \mu_2+\mu_1
q_{21}+\mu_2 q_{12}}.\label{ES2}
\end{eqnarray}

\subsection{The fully unobservable case}

We can now proceed and determine the equilibrium strategies of the
customers in the fully unobservable case. A general balking strategy
in the fully unobservable case is specified by a single joining
probability $q$. The case $q=0$ corresponds to the pure strategy `to
balk' whereas the case $q=1$ corresponds to the pure strategy `to
join'. Any value of $q\in (0,1)$ corresponds to a mixed (randomized)
strategy `to join with probability $q$ or balk with probability
$1-q$'. We have the following Theorem
\ref{equilibrium-fully-unobservable-theorem}.

\begin{thm}\label{equilibrium-fully-unobservable-theorem}
In the fully unobservable model of the stochastic clearing system in
alternating environment, there always exists a dominant strategy. The dominant strategies depend on the relative value of the ratio $\frac{R}{C}$ with respect to the critical value
\begin{equation}
V_{fu}=\frac{\lambda_1 q_{21} \mu_2 +\lambda_2 q_{12}
\mu_1}{(\lambda_1 q_{21} + \lambda_2 q_{12})(\mu_1 \mu_2 +\mu_1
q_{21} + \mu_2 q_{12})}+\frac{q_{21}+q_{12}}{\mu_1 \mu_2 + \mu_1
q_{21} + \mu_2 q_{12}}.
\end{equation}
We have three cases that are summarized in Table \ref{Table-Unobservable-Strategies}.
\medskip
\begin{table}[h]
\begin{center}
\begin{tabular}{|l||c|c|c|}
\hline
& & & \\
Value of $\frac{R}{C}$ & $\frac{R}{C}<V_{fu}$ & $\frac{R}{C}=V_{fu}$ & $\frac{R}{C}>V_{fu}$\\
& & & \\
\hline
& & & \\
Dominant strategy(ies) & $0$ & $q\in [0,1]$ & $1$\\
& & & \\
\hline
\end{tabular} \caption{Dominant strategies in the fully unobservable case}\label{Table-Unobservable-Strategies}
\end{center}
\end{table}
\end{thm}
\vspace{-1cm}

\begin{prooof}{Proof} Suppose that the customers follow
a certain strategy and consider a tagged customer upon arrival. The
probability that he finds the environment at state $e$ is
\begin{equation}
p_{E-arrival}(e)=\frac{\lambda_e p_E(e)}{
(\lambda_1 p_E(1)+\lambda_2 p_E(2))},\label{p-E-arrival}
\end{equation}
where $(p_E(e),e=1,2)$ is the stationary distribution of the environment which is given by
\begin{eqnarray}
p_E(1)&=&\frac{q_{21}}{q_{21}+q_{12}}\label{env-stat-1}\\
p_E(2)&=&\frac{q_{12}}{q_{21}+q_{12}}\label{env-stat-2}.
\end{eqnarray}
Therefore, the expected net benefit of
the tagged customer if he decides to join is given by
\begin{equation}
S_{fu}=\sum_{e=1}^2 p_{E-arrival}(e)
(R-CE[S_{e}]),\label{S-un-1}
\end{equation}
where $E[S_{e}]$ are given by \eqref{ES1}-\eqref{ES2}. Plugging \eqref{env-stat-1}-\eqref{env-stat-2} in \eqref{p-E-arrival} and substituting in \eqref{S-un-1} yields
\begin{eqnarray}
S_{fu}&=&R-C\;\frac{\lambda_1 q_{21} E[S_{1}]+\lambda_2 q_{12} E[S_{2}]}{\lambda_1 q_{21} + \lambda_2 q_{12}}\nonumber\\
&=&R-C \left( \frac{\lambda_1 q_{21} \mu_2 +\lambda_2 q_{12}
\mu_1}{(\lambda_1 q_{21} + \lambda_2 q_{12})(\mu_1 \mu_2 +\mu_1
q_{21} + \mu_2 q_{12})}+\frac{q_{21}+q_{12}}{\mu_1 \mu_2 + \mu_1
q_{21} + \mu_2 q_{12}} \right).
\end{eqnarray}
The tagged customer prefers to join if $S_{fu}> 0$, prefers to balk
if $S_{fu}<0$ and he is indifferent between joining and balking if
$S_{fu}=0$ . Solving with respect to $\frac{R}{C}$, we obtain the
three cases of Table \ref{Table-Unobservable-Strategies}.
\end{prooof}\\

\bigskip
\bigskip
\bigskip

\subsection{The almost unobservable case}

We can now proceed and determine the equilibrium strategies of the customers in the almost unobservable case. A general balking strategy in the almost unobservable case is specified by an ordered pair of joining probabilities $(q_1,q_2)$, where $q_e$ is the joining probability of a customer if the environmental state upon arrival is $e$, $e=1,2$. We have the following Theorem \ref{equilibrium-almost-unobservable-theorem}.

\bigskip

\begin{thm}\label{equilibrium-almost-unobservable-theorem}
In the almost unobservable model of the stochastic clearing system in alternating environment, there always exists a dominant strategy. The dominant strategies depend on the relative value of the ratio $\frac{R}{C}$ with respect to the critical values
\begin{eqnarray}
V_{au}^{min}=\frac{\min(\mu_1,\mu_2)+q_{21}+q_{12}}{\mu_1 \mu_2+\mu_1 q_{21}+\mu_2 q_{12}},\;
V_{au}^{max}=\frac{\max(\mu_1,\mu_2)+q_{21}+q_{12}}{\mu_1 \mu_2+\mu_1 q_{21}+\mu_2 q_{12}}.
\end{eqnarray}

If $\mu_1\ne \mu_2$, then $V_{au}^{min}<V_{au}^{max}$ and we have five cases that are summarized in Table \ref{Table-Almost-Unobservable-Strategies-1}.

\medskip
\begin{table}[h]
\begin{center}
\begin{tabular}{|l||c|c|c|c|c|}
\hline
& & & & &\\
Value of $\frac{R}{C}$ & $\frac{R}{C}<V_{au}^{min}$ & $\frac{R}{C}=V_{au}^{min}$ & $V_{au}^{min}<\frac{R}{C}<V_{au}^{max}$ & $\frac{R}{C}=V_{au}^{max}$ & $\frac{R}{C}>V_{au}^{max}$\\
\hline
& & & & &\\
Dominant & & & & &  \\
strategy(ies),& $(0,0)$ & $(0,q_2),$ & $(0,1)$ & $(q_1,1),$ & $(1,1)$\\
when $\mu_1<\mu_2$ & & $q_2\in [0,1]$ & & $q_1\in [0,1]$ &\\
& & & & &\\
\hline
& & & & &\\
Dominant & & & & &  \\
strategy(ies),& $(0,0)$ & $(q_1,0),$ & $(1,0)$ & $(1,q_2),$ & $(1,1)$\\
when $\mu_1>\mu_2$ & & $q_1\in [0,1]$ & & $q_2\in [0,1]$ &\\
& & & & &\\
\hline
\end{tabular} \caption{Dominant strategies in the almost unobservable case, when $\mu_1\ne \mu_2$}\label{Table-Almost-Unobservable-Strategies-1}
\end{center}
\end{table}
\vspace{-0.5cm}
\newpage

If $\mu_1= \mu_2$, then $V_{au}^{min}=V_{au}^{max}$. Let $V_{au}$ denote the common value of $V_{au}^{min}$ and $V_{au}^{max}$. We have three cases that are summarized in Table \ref{Table-Almost-Unobservable-Strategies-2}.
\medskip
\begin{table}[h]
\begin{center}
\begin{tabular}{|l||c|c|c|}
\hline
& & & \\
Value of $\frac{R}{C}$ & $\frac{R}{C}<V_{au}$ & $\frac{R}{C}=V_{au}$ & $\frac{R}{C}>V_{au}$\\
& & & \\
\hline
& & & \\
Dominant strategy(ies), & $(0,0)$ & $(q_1,q_2),$ & $(1,1)$\\
when $\mu_1= \mu_2$ & & $q_1,q_2\in [0,1]$ &\\
& & & \\
\hline
\end{tabular} \caption{Dominant strategies in the almost unobservable case, when $\mu_1 = \mu_2$ }\label{Table-Almost-Unobservable-Strategies-2}
\end{center}
\end{table}
\end{thm}
\vspace{-1cm}
\begin{prooof}{Proof} Consider a tagged customer that observes the
state of the environment upon arrival. If he decides to join given
that he finds the environment at state $e$, then his expected net
benefit will be
\begin{equation}
S_{au}(e)=R-CE[S_{e}],\label{S-obs-ne}
\end{equation}
where $E[S_{e}]$ are given by
\eqref{ES1}-\eqref{ES2}. The customer prefers to join if $S_{au}(e)>
0$, which is written equivalently as $\frac{R}{C}> E[S_{e}]$.
Similarly, we have that he prefers to balk if $\frac{R}{C}<
E[S_{e}]$ and he is indifferent between joining and balking if
$\frac{R}{C}= E[S_{e}]$. By considering the various possible cases
with regard to the order of the three quantities $\frac{R}{C}$,
$E[S_{1}]$ and $E[S_{2}]$, we obtain the corresponding cases in
the statement of the Theorem \ref{equilibrium-almost-unobservable-theorem}. Note that the strategies prescribed in
the Theorem \ref{equilibrium-almost-unobservable-theorem} are dominant, since they do not depend on what the other
customers do, i.e. they are best responses against any strategy of
the others.
\end{prooof}

\subsection{The fully observable case}

Regarding the fully observable case, where the arriving customers observe both the number of waiting customers and the state of the environment, the situation is identical to the almost unobservable case. This happens because the
mean sojourn time of an arriving customer, given that he finds $n$ customers in the system and the environment at state $e$ does not depend on $n$. Therefore, if the environmental state $e$ is observed upon arrival, then the information about the number of customers $n$ is superfluous and is discarded by the customers. We conclude that the dominant balking strategies
are the ones described in Theorem \ref{equilibrium-almost-unobservable-theorem}.\\

\section{The almost observable case: Preliminaries}\label{Observable-Section-Preliminaries}

In this section, we consider the almost observable case. In this case, the customers, upon arrival and before making their decisions about whether to join or balk, observe the number of customers in the system but not the state of the environment. Thus a general balking strategy in this case is specified by a vector of joining probabilities $(\theta_0,\theta_1,\theta_2,\ldots)$, where
$\theta_i$ is the joining probability of a customer that sees $i$
customers in the system upon arrival (excluding himself).\\

Suppose that a tagged customer observes $n$ customers in the system upon arrival. Although his mean sojourn time does not depend on $n$, the information about $n$ influences the probabilities that the environment is found at state 1 or 2. We expect intuitively that there are two cases: Either the `slow service' environmental state $e$ with $\mu_e=\min(\mu_1,\mu_2)$
coincides to the `more congested' environmental state $e'$ with
$\rho_{e'}=\max
(\rho_1,\rho_2)$ or it coincides
to the `less congested' environmental state $e''$ with
$\rho_{e''}=\min
(\rho_1,\rho_2)$.
In the former case, the greater the number $n$ of the customers
found by the tagged customer, the more probable is that the
environment is found at the `slow service' environmental state.
Therefore, the tagged customer becomes less willing to join the
system as $n$ increases. Thus, we expect that the tagged customer
will benefit from joining the system, if the number of customers
$n$ is below a certain threshold, i.e. he will adopt a threshold
strategy. On the contrary, in the latter case, the situation is
reversed. Then, the greater the number $n$ of the customers found by
a tagged customer, the more probable is that the environment is
found at the `fast service' environmental state. Therefore, we
expect that the tagged customer will benefit from joining the
system, if the number of customers $n$ exceeds a certain threshold,
i.e. he will adopt a so called reverse-threshold strategy. Following
this reasoning, we will limit our search for equilibrium strategies
within the class of threshold and reverse-threshold strategies. As
we will see, this family is rich enough to ensure the existence of
an equilibrium strategy for any values of the underlying parameters
of the model.

\begin{defn}\label{defin-strat} A balking strategy $(\theta_0,\theta_1,\theta_2,\ldots)$, where
$\theta_i$ is the joining probability of a customer that sees $i$
customers in the system upon arrival (excluding himself) is said to
be a mixed threshold strategy, if there exist $n_0\in\{0,1,\ldots\}$
and $\theta\in[0,1]$ such that $\theta_i=1$, for $i<n_0$,
$\theta_{n_0}=\theta$ and $\theta_i=0$, for $i>n_0$. Such a strategy
will be referred to as the $(n_0,\theta)$-mixed threshold strategy
(symbolically the $\lceil n_0,\theta \rceil$ strategy) and it
prescribes to join if you see less than $n_0$ customers, to join
with probability $\theta$ if you see exactly $n_0$ customers and to
balk if you see more than $n_0$ customers.

An $(n_0,0)$-mixed threshold strategy which prescribes to join if
you see less than $n_0$ customers and to balk otherwise will be
referred to as the $n_0$-pure threshold strategy (symbolically the
$\lceil n_0 \rceil$ strategy).

A balking strategy $(\theta_0,\theta_1,\theta_2,\ldots)$ is said to
be a mixed reverse-threshold strategy, if there exist
$n_0\in\{0,1,\ldots\}$ and $\theta\in[0,1]$ such that $\theta_i=0$,
for $i<n_0$, $\theta_{n_0}=\theta$ and $\theta_i=1$, for $i>n_0$.
Such a strategy will be referred to as the $(n_0,\theta)$-mixed
reverse-threshold strategy (symbolically the
$\lfloor n_0,\theta \rfloor $ strategy) and it prescribes to balk if you
see less than $n_0$ customers, to join with probability $\theta$ if
you see exactly $n_0$ customers and to join if you see more that
$n_0$ customers.

An $(n_0,1)$-mixed reverse-threshold strategy which prescribes to
join if you see at least $n_0$ customers and to balk otherwise will
be referred to as the $n_0$-pure reverse-threshold strategy
(symbolically the $\lfloor n_0 \rfloor $ strategy).

The strategy which prescribes to join in any case is considered to
be both a threshold and a reverse-threshold strategy (symbolically
the $\lceil \infty \rceil$ or $\lfloor 0 \rfloor $ strategy). The same is
true for the strategy which prescribes to balk in any case
(symbolically the $\lceil 0 \rceil$ or $\lfloor \infty \rfloor $ strategy).
\end{defn}

\subsection{Stationary distributions}\label{stat-distr-subsec}

In this subsection, we determine the stationary distributions of the system, when the customers follow any given strategy from the ones that have been described in Definition \ref{defin-strat}. We will first determine the stationary distribution of the original system when all customers join. The result is reported in the following Proposition \ref{stationary-original}.

\begin{prop}\label{stationary-original}
Consider the stochastic clearing system in alternating environment,
where all customers join. The stationary distribution $(p(n,e))$ is
given by the formulas
\begin{eqnarray}
p(n,1)&=&A_1\left(\frac{1}{1-z_{1}}\right)^{n}+B_1\left(\frac{1}{1-z_{2}}\right)^{n},\ \ \ n\geq 0,\label{stat-prob-1}\\
p(n,2)&=&A_2\left(\frac{1}{1-z_{1}}\right)^{n}+B_2\left(\frac{1}{1-z_{2}}\right)^{n},\label{stat-prob-2}\
\ \ n\geq 0,
\end{eqnarray}
where
\begin{eqnarray}
A_1&=&\frac{(\mu_{1}\lambda_{2}z_{1}+\mu_{1}\mu_{2}+\mu_{2}q_{12}+\mu_{1}q_{21})p_E(1)}{\sqrt{\Delta}(1-z_1)},\label{A1}\\
B_1&=&-\frac{(\mu_{1}\lambda_{2}z_{2}+\mu_{1}\mu_{2}+\mu_{2}q_{12}+\mu_{1}q_{21})p_E(1)}{\sqrt{\Delta}(1-z_2)},\label{B1}\\
A_2&=&\frac{(\mu_{2}\lambda_{1}z_{1}+\mu_{1}\mu_{2}+\mu_{2}q_{12}+\mu_{1}q_{21})p_E(2)}{\sqrt{\Delta}(1-z_1)},\label{A2}\\
B_2&=&-\frac{(\mu_{2}\lambda_{1}z_{2}+\mu_{1}\mu_{2}+\mu_{2}q_{12}+\mu_{1}q_{21})p_E(2)}{\sqrt{\Delta}(1-z_2)},\label{B2}\\
\Delta&=&[\lambda_2(\mu_1+q_{12})-\lambda_1(\mu_2+q_{21})]^2+4\lambda_1\lambda_2q_{12}q_{21},\label{Delta}\\
z_{1,2}&=&\frac{-\lambda_1(\mu_2+q_{21})-\lambda_2(\mu_1+q_{12})\pm\sqrt{\Delta}}{2\lambda_1\lambda_2}\label{z12}
\end{eqnarray}
and $p_E(1),\ p_E(2)$ are the stationary probabilities of $\{E(t)\}$
given from \eqref{env-stat-1}-\eqref{env-stat-2}.
\end{prop}
\begin{prooof}{Proof}
For the stationary analysis, note that the state of the system is
described by a continuous-time Markov chain with state space
$S^{N,E}=\{(n,e): n\geq 0,\ e=\ 1,\ 2\}$ with its non-zero
transition rates given by \eqref{TranRate1}-\eqref{TranRate4}. The
corresponding stationary distribution $(p(n,e):(n,e)\in S^{N,E})$ is
obtained as the unique positive normalized solution of the following
system of balance equations:
\begin{eqnarray}
(\lambda_1+\mu_{1}+q_{12})p(0,1)&=&q_{21}p(0,2)+\sum_{n=0}^{\infty}\mu_1p(n,1),\label{be11}\\
(\lambda_1+\mu_{1}+q_{12})p(n,1)&=&q_{21}p(n,2)+\lambda_1p(n-1,1),\ \  n\geq 1,\label{be21}\\
(\lambda_2+\mu_{2}+q_{21})p(0,2)&=&q_{12}p(0,1)+\sum_{n=0}^{\infty}\mu_2p(n,2),\label{be12}\\
(\lambda_2+\mu_{2}+q_{21})p(n,2)&=&q_{12}p(n,1)+\lambda_2p(n-1,2),\
\   n\geq 1,\label{be22}
\end{eqnarray}
where we have included in \eqref{be11}, \eqref{be12} the
pseudo-transitions from $(0,e)$ to $(0,e),\ e=1, 2$, with rate
$\mu_{e}$, that correspond to visits of the facility at an empty
system. Note also that the underlying Markov chain is always
positive recurrent as the stochastic clearing mechanism ensures that
starting from the state $(0,1)$, the process will visit it again
with
probability 1 and the corresponding mean time is finite.\\

For determining the stationary probabilities, we may follow the standard probability generating function approach. Thus, we define the partial stationary probability generating functions of the system as
\begin{equation}
G_{e}(z)=\sum_{n=0}^{\infty}p(n,e)z^n,\ \ |z|\leq 1,\ \ e=1,2.\label{gf-original}
\end{equation}
and we observe that $G_1(1)=p_{E}(1)$, $G_2(1)=p_{E}(2)$ with $p_{E}(1)$,
$p_{E}(2)$ given from \eqref{env-stat-1}-\eqref{env-stat-2}. Summing
equation \eqref{be11} and equations \eqref{be21} multiplied by
$z^n$, $n\geq 1$, yields after some straightforward algebra a linear equation in $G_{1}(z)$ and $G_{2}(z)$. Similarly, equations \eqref{be12} and \eqref{be22}, $n\geq 1$, yield another linear equation in $G_{1}(z)$ and $G_{2}(z)$. Solving the system of these equations we obtain $G_{1}(z)$ and $G_{2}(z)$ as rational functions of $z$ with known coefficients expressed in terms of the parameters of the model. Using partial fraction expansion and then expanding the simple fractions in powers of $z$ yields \eqref{stat-prob-1} and \eqref{stat-prob-2}. Indeed, by direct substitution, we can easily check that $p(n,1)$ and $p(n,2)$ given by \eqref{stat-prob-1} and \eqref{stat-prob-2} satisfy \eqref{be21}. By a simple summation, we can also check that $p(n,1)$ and $p(n,2)$ given by \eqref{stat-prob-1} and \eqref{stat-prob-2} satisfy \eqref{be11}. The validity of \eqref{be22} and \eqref{be12} is checked similarly.
\end{prooof}\\

We will now deduce the stationary distribution of the system when
the customers follow a mixed threshold strategy. We have the
following Proposition \ref{stationary-mixed}.

\begin{prop}\label{stationary-mixed}
Consider the almost observable model of the stochastic clearing
system in alternating environment, where the customers join the
system according to the $(n_{0},\theta)$-mixed threshold strategy.
The corresponding stationary distribution
$(p_{ao}(n,e;\lceil n_{0},\theta \rceil))$ is given by the formulas
\begin{eqnarray}
p_{ao}(n,e;\lceil n_{0},\theta \rceil)&=&p(n,e),\ \ \ 0 \leq n \leq n_{0}-1, \ \ \ e=1,2, \label{stat-prob-mix-n0-}\\
p_{ao}(n_{0},e;\lceil n_{0},\theta \rceil)&=&\sum_{n=n_{0}}^{\infty}(1-\theta)^{n-n_{0}}p(n,e),\ \ e=1,2,  \label{stat-prob-mix-n0}\\
p_{ao}(n_{0}+1,e;\lceil n_{0},\theta \rceil)&=&\sum_{n=n_{0}+1}^{\infty}[1-(1-\theta)^{n-n_{0}}]p(n,e),\ \ e= 1,2, \label{stat-prob-mix-n0+1} \\
p_{ao}(n,e;\lceil n_{0},\theta \rceil)&=&0,\ \ \ n \geq n_{0}+2,\ \
e=1,2,\label{stat-prob-mix-n0+1+}
\end{eqnarray}
where $p(n,e)$ are given by \eqref{stat-prob-1}-\eqref{stat-prob-2}.
\end{prop}

\begin{prooof}{Proof}
We assume that the customers follow the $(n_{0},\theta)$-mixed
threshold strategy. Then the evolution of the system can be
described by a Markov chain which is absorbed with probability 1 in
the positive recurrent closed class of states
$S_{ao}^{N,E}(\lceil n_{0},\theta \rceil)=\{(n,e):\; 0\leq n\leq n_{0}
+1,\ e=\ 1,\ 2\}$. For the sake of brevity, we suppress the notation
regarding $\lceil n_{0},\theta \rceil$ in the rest of the proof . Thus,
we will refer to the corresponding stationary probabilities
$p_{ao}(n,e;\lceil n_{0},\theta \rceil)$ by $p_{ao}(n,e)$.\\

Since the Markov chain is finally absorbed in
$S_{ao}^{N,E}(\lceil n_{0},\theta \rceil)$, we obtain immediately
\eqref{stat-prob-mix-n0+1+}. The vector of the stationary
probabilities $(p_{ao}(n,e):(n,e)\in
S_{ao}^{N,E}(\lceil n_{0},\theta \rceil))$ is obtained as the unique
positive normalized solution of the system of balance equations
\begin{eqnarray}
(\lambda_1+\mu_{1}+q_{12})p_{ao}(0,1)&=&q_{21}p_{ao}(0,2)+\sum_{n=0}^{n_0+1}\mu_1p_{ao}(n,1),\label{bem11}\\
(\lambda_1+\mu_{1}+q_{12})p_{ao}(n,1)&=&q_{21}p_{ao}(n,2)+\lambda_1p_{ao}(n-1,1),\ 1 \leq n \leq n_{0}-1,\label{bem21}\\
(\lambda_1\theta+\mu_{1}+q_{12})p_{ao}(n_0,1)&=&q_{21}p_{ao}(n_0,2)+\lambda_1p_{ao}(n_0-1,1),\label{bem31}\\
(\mu_{1}+q_{12})p_{ao}(n_0+1,1)&=&q_{21}p_{ao}(n_0+1,2)+\lambda_1\theta  p_{ao}(n_0,1),\label{bem41}\\
(\lambda_2+\mu_{2}+q_{21})p_{ao}(0,2)&=&q_{12}p_{ao}(0,1)+\sum_{n=0}^{n_0+1}\mu_2p_{ao}(n,2),\label{bem12}\\
(\lambda_2+\mu_{2}+q_{21})p_{ao}(n,2)&=&q_{12}p_{ao}(n,1)+\lambda_2p_{ao}(n-1,2),\ 1 \leq n \leq n_{0}-1,\label{bem22}\\
(\lambda_2\theta+\mu_{2}+q_{21})p_{ao}(n_0,2)&=&q_{12}p_{ao}(n_0,1)+\lambda_2p_{ao}(n_0-1,2),\label{bem32}\\
(\mu_{2}+q_{21})p_{ao}(n_0+1,2)&=&q_{12}p_{ao}(n_0+1,1)+\lambda_2\theta
p_{ao}(n_0,2),\label{bem42}
\end{eqnarray}
where we have included in \eqref{bem11} and \eqref{bem12} the
pseudo-transitions from $(0,e)$ to $(0,e),\ e=1,2$, with rate
$\mu_{e}$, that correspond to visits of the facility at an empty system.

For deducing the formulas \eqref{stat-prob-mix-n0-}-\eqref{stat-prob-mix-n0+1+} for the stationary probabilities, we may again follow the standard probability generating function approach, as it was briefly described in the proof of Proposition \ref{stationary-original}. However, given the formulas \eqref{stat-prob-mix-n0-}-\eqref{stat-prob-mix-n0+1+}, it is easy to check by direct substitution that the stationary probabilities satisfy the equations \eqref{bem11}-\eqref{bem42} (using also simple summations for \eqref{bem11} and \eqref{bem12}).
\end{prooof}\\

We can now readily conclude the following Corollaries
\ref{stationary-pure} and \ref{stationary-pure-0}.

\begin{cor}\label{stationary-pure}
Consider the almost observable model of the stochastic clearing
system in alternating environment, where the customers join the
system according to the $n_{0}$-pure threshold strategy. The
corresponding stationary distribution
$(p_{ao}(n,e;\lceil n_{0} \rceil ))$ is given by the formulas
\begin{eqnarray}
p_{ao}(n,e;\lceil n_{0} \rceil )&=&p(n,e),\ \ \ 0 \leq n \leq n_{0}-1, \ \ \ e=1,2, \label{stat-prob-pure-n0-}\\
p_{ao}(n_{0},e;\lceil n_{0} \rceil )&=&\sum_{n=n_{0}}^{\infty}p(n,e),\ \ e=1,2,  \label{stat-prob-pure-n0}\\
p_{ao}(n,e;\lceil n_{0} \rceil )&=&0,\ \ \ n \geq n_{0}+1,\ \
e=1,2,\label{stat-prob-pure-n0+}
\end{eqnarray}
where $p(n,e)$ are given by \eqref{stat-prob-1}-\eqref{stat-prob-2}.
\end{cor}

\begin{cor}\label{stationary-pure-0}
Consider the almost observable model of the stochastic clearing
system in alternating environment, where the customers always balk.
The corresponding stationary distribution
$(p_{ao}(n,e;\lceil 0 \rceil ))$ is given by the formulas
\begin{eqnarray}
p_{ao}(0,e;\lceil 0 \rceil )&=&p_E(e),\ \ e=1,2,  \label{stat-prob-pure-0}\\
p_{ao}(n,e;\lceil 0 \rceil )&=&0,\ \ \ n \geq 1,\ \
e=1,2,\label{stat-prob-pure-0+}
\end{eqnarray}
where $p_E(e)$, $e=1,2$ are given by \eqref{env-stat-1}-\eqref{env-stat-2}.
\end{cor}

We will now deduce the stationary distribution of the system when
the customers follow an $(n_0,\theta)$-mixed reverse-threshold
strategy.

\begin{remark}\label{stationary-reverse-trivial}
Under an $(n_0,\theta)$-mixed reverse-threshold strategy with
$n_0\geq 1$, we have that the customers balk when they arrive at an
empty system. Thus we have the stationary distribution of Corollary
\ref{stationary-pure-0}.
\end{remark}

It is left to show what happens when the customers follow a $(0,\theta)$-mixed reverse-threshold strategy.

\begin{prop}\label{stationary-reverse-main}
Consider the almost observable model of the stochastic clearing
system in alternating environment, where the customers join the
system according to a $(0,\theta)$-mixed reverse-threshold strategy.
For $\theta=0$, the stationary distribution
$(p_{ao}(n,e;\lfloor 0,0 \rfloor ))$ is given by the formulas
\begin{eqnarray}
p_{ao}(0,e;\lfloor 0,0 \rfloor )&=&p_E(e),\ \ e=1,2,  \label{stat-rev-prob-pure-0}\\
p_{ao}(n,e;\lfloor 0,0 \rfloor )&=&0,\ \ \ n \geq 1,\ \
e=1,2,\label{stat-rev-prob-pure-0+}
\end{eqnarray}
where $p_E(e)$, $e=1,2$ are given by \eqref{env-stat-1}-\eqref{env-stat-2}.

For $\theta\in(0,1)$, the stationary distribution
$(p_{ao}(n,e;\lfloor 0,\theta \rfloor ))$ is given by the formulas
\begin{eqnarray}
p_{ao}(0,e;\lfloor 0,\theta \rfloor )&=&\sum_{n=0}^{\infty}(1-\theta)^np(n,e),\
\ e=1,2,\label{stat-rev-prob-main-0}\\
p_{ao}(n,e;\lfloor 0,\theta \rfloor )&=&\theta \sum_{i=n}^{\infty}
(1-\theta)^{i-n} p(i,e),\ \ n\geq 1,\ \ e=1,2,\label{stat-rev-prob-main-0+}
\end{eqnarray}
where $p(n,e)$ are given by \eqref{stat-prob-1}-\eqref{stat-prob-2}.

For $\theta=1$, the stationary distribution
$(p_{ao}(n,e;\lfloor 0,1 \rfloor ))$ is given by the formula
\begin{equation}
p_{ao}(n,e;\lfloor 0,1 \rfloor )=p(n,e),\ \ n\geq 0,\ \ e=1,2,
\label{stat-rev-prob-always-join}
\end{equation}
where $p(n,e)$ are given by \eqref{stat-prob-1}-\eqref{stat-prob-2}.
\end{prop}

The proof of Proposition \ref{stationary-reverse-main} for
$\theta=0$ is immediate, as in this case the customers balk whenever
they arrive at an empty system. Therefore under such a strategy the
corresponding continuous-time Markov chain is absorbed with
probability 1 into the subset $\{(0,1),(0,2)\}$ of the state space
and the stationary distribution is the one given by
\eqref{stat-prob-pure-0} and \eqref{stat-prob-pure-0+} as in
Corollary \ref{stationary-pure-0}. In case $\theta=1$, the customers
always join so we apply Proposition \ref{stationary-original}. Thus,
the only interesting case is for $\theta\in(0,1)$. Then, the proof
of Proposition \ref{stationary-reverse-main} follows a similar line
of argument as the proofs of Propositions \ref{stationary-original}
and \ref{stationary-mixed}. Therefore, for the sake of brevity, it
is omitted.

\subsection{Expected net benefit functions}

Based on the results of subsection
\ref{stat-distr-subsec}, we can now compute the
expected net benefit of a tagged customer if he decides to join the
system after observing $n$ customers upon arrival. Of course, his
expected net benefit depends on the strategy followed by the other
customers. Thus, we have various cases, according to whether the
customers follow a threshold or a reverse-threshold strategy. We
have the following Propositions \ref{expected-benefit-all-enter}--\ref{expected-benefit-reverse} and the Corollary \ref{expected-benefit-pure}.

\begin{prop}\label{expected-benefit-all-enter}
Consider the almost observable model of the stochastic clearing
system in alternating environment, where all customers join the
system.  Then, the expected net benefit
$S_{ao}(n;\lceil \infty \rceil )\equiv S_{ao}(n;\lfloor 0 \rfloor )$ of an arriving customer, if he decides to join, given that he
finds $n$ customers in the system, is given by
\begin{eqnarray}
S_{ao}(n;\lceil \infty \rceil )\equiv
S_{ao}(n;\lfloor 0 \rfloor )&=&R-C\frac{A\left(\frac{1}{1-z_1}\right)^n+B\left(\frac{1}{1-z_2}\right)^n}{D\left(\frac{1}{1-z_1}\right)^n+E\left(\frac{1}{1-z_2}\right)^n},\
\ \ \ n\geq 0,\label{exp-ben}
\end{eqnarray}
where
\begin{eqnarray}
A&=&\lambda_1 A_1E[S_{1}]+\lambda_2 A_2E[S_{2}],\label{A}\\
B&=&\lambda_1 B_1E[S_{1}]+\lambda_2 B_2E[S_{2}],\label{B}\\
D&=&\lambda_1 A_1+\lambda_2 A_2,\label{D}\\
E&=&\lambda_1B_1+\lambda_2 B_2\label{E}
\end{eqnarray}
and $E[S_{1}],\ E[S_{2}],\  A_1,\ B_1,\ A_2,\  B_2,\ z_1,\ z_2$
are given by \eqref{ES1}-\eqref{ES2}, \eqref{A1}-\eqref{B2} and
\eqref{z12}.
\end{prop}

\begin{prooof}{Proof}
The mean sojourn time of an arriving customer, if he decides
to join, given that he finds $n$ customers in the system is given by
\begin{equation}
p^{-}_{ao}(1|n;\lceil \infty \rceil )E[S_{1}]+p^{-}_{ao}(2|n;\lceil \infty \rceil )E[S_{2}],
\end{equation}
where $p^{-}_{ao}(e|n;\lceil \infty \rceil )$, $e=1,2,$ is the
probability that an arriving customer finds the environment at state
$e$, given that he observes $n$ customers in the system and that the
$\lceil \infty \rceil $-strategy is followed by the other customers. The embedded (Palm)
probabilities $p^{-}_{ao}(e|n;\lceil \infty \rceil )$ are given by
\begin{equation}\label{p-ao-unrestricted}
p^{-}_{ao}(e|n;\lceil \infty \rceil )=\frac{\lambda_e
p(n,e)}{\lambda_1p(n,1)+\lambda_2p(n,2)},\ e=1,2,
\end{equation}
where $p(n,e)$ are given by \eqref{stat-prob-1}-\eqref{stat-prob-2}.
Thus, the expected benefit of the tagged arriving customer, if he
decides to join, is equal to
\begin{equation}
S_{ao}(n;\lceil \infty \rceil )\equiv
S_{ao}(n;\lfloor 0 \rfloor )=R-C[p^{-}_{ao}(1|n;\lceil \infty \rceil )E[S_{1}]+p^{-}_{ao}(2|n;\lceil \infty \rceil )E[S_{2}]].\label{Sao-unrestricted}
\end{equation}
Plugging the formulas \eqref{stat-prob-1}-\eqref{stat-prob-2} into \eqref{p-ao-unrestricted} and subsequently into
\eqref{Sao-unrestricted} yields \eqref{exp-ben}.
\end{prooof}\\

\begin{prop}\label{expected-benefit-mixed}
Consider the almost observable model of the stochastic clearing
system in alternating environment, where the customers join the
system according to the $(n_{0},\theta)$-mixed threshold strategy.
Then, the expected net benefit $S_{ao}(n;\lceil n_{0},\theta \rceil )$
of an arriving customer, if he decides to join, given that he
finds $n$ customers in the system, is given by
\begin{eqnarray}
S_{ao}(n;\lceil n_{0},\theta \rceil )&=&R-C\frac{A\left(\frac{1}{1-z_1}\right)^n+B\left(\frac{1}{1-z_2}\right)^n}{D\left(\frac{1}{1-z_1}\right)^n+E\left(\frac{1}{1-z_2}\right)^n},\ \ \ 0 \leq n \leq n_{0}-1,\label{exp-ben-mix-n0-}\\
S_{ao}(n_0;\lceil n_{0},\theta \rceil )&=&R-C\frac{\sum_{i=n_0}^{\infty}(1-\theta)^{i-n_0}\left[A\left(\frac{1}{1-z_1}\right)^i+B\left(\frac{1}{1-z_2}\right)^i\right]}{\sum_{i=n_0}^{\infty}(1-\theta)^{i-n_0}\left[D\left(\frac{1}{1-z_1}\right)^i+E\left(\frac{1}{1-z_2}\right)^i\right]},
\label{exp-ben-mix-n0}\\
S_{ao}(n_0+1;\lceil n_{0},\theta \rceil )&=&R-C\frac{\sum_{i=n_0+1}^{\infty}[1-(1-\theta)^{i-n_0}]\left[A\left(\frac{1}{1-z_1}\right)^i+B\left(\frac{1}{1-z_2}\right)^i\right]}{\sum_{i=n_0+1}^{\infty}[1-(1-\theta)^{i-n_0}]\left[D\left(\frac{1}{1-z_1}\right)^i+E\left(\frac{1}{1-z_2}\right)^i\right]},
\label{exp-ben-mix-n0+1}
\end{eqnarray}
where
 $A,\ B,\ D,\ E,\ z_1,\ z_2$ are given by \eqref{A}-\eqref{E} and
\eqref{z12}.
\end{prop}

\begin{prooof}{Proof}
Assume that the customers join the system according to the
$(n_{0},\theta)$-mixed threshold strategy. Then, the mean sojourn
time of a tagged arriving customer, if he decides to join, given
that he finds $n$ customers in the system is given by
\begin{equation}
p^{-}_{ao}(1|n;\lceil n_{0},\theta \rceil )E[S_{1}]+p^{-}_{ao}(2|n;\lceil n_{0},\theta \rceil )E[S_{2}],
\end{equation}
where $p^{-}_{ao}(e|n;\lceil n_{0},\theta \rceil ),\ e=1,2,$ is the
probability that an arriving customer finds the environment at state
$e$, given that there are $n$ customers in the system and that the
$\lceil n_0,\theta \rceil $-strategy is followed. The embedded (Palm)
probabilities are given by
\begin{equation}
p^{-}_{ao}(e|n;\lceil n_0,\theta \rceil )=\frac{\lambda_e
p_{ao}(n,e;\lceil n_{0},\theta \rceil )}{\lambda_1p_{ao}(n,1;\lceil n_0,\theta \rceil )+\lambda_2p_{ao}(n,2;\lceil n_0,\theta \rceil )},\
e=1,2,
\end{equation}
where $p_{ao}(n,e;\lceil n_0,\theta \rceil )$ are given by
\eqref{stat-prob-mix-n0-}-\eqref{stat-prob-mix-n0+1}. Thus, the
expected benefit of the tagged customer, if he decides to join, is
equal to
\begin{equation}
S_{ao}(n;\lceil n_{0},\theta \rceil )=R-C[p^{-}_{ao}(1|n;\lceil n_{0},\theta \rceil )E[S_{1}]+p^{-}_{ao}(2|n;\lceil n_{0},\theta \rceil )E[S_{2}]].
\end{equation}
Using the various forms of $p_{ao}(n,e;\lceil n_{0},\theta \rceil )$ in
\eqref{stat-prob-mix-n0-}-\eqref{stat-prob-mix-n0+1} yields \eqref{exp-ben-mix-n0-}-\eqref{exp-ben-mix-n0+1}.
\end{prooof}\\

In the case of the $n_{0}$-pure threshold strategy, we obtain the
following Corollary \ref{expected-benefit-pure}.

\begin{cor}\label{expected-benefit-pure}
Consider the almost observable model of the stochastic clearing
system in alternating environment, where the customers join the
system according to the $n_0$-pure threshold strategy. Then, the
expected net benefit $S_{ao}(n;\lceil n_0 \rceil )$ of an arriving
customer, if he decides to join, given that he finds $n$ customers
in the system, is given by
\begin{eqnarray}
S_{ao}(n;\lceil n_0 \rceil )&=&R-C\frac{A\left(\frac{1}{1-z_1}\right)^n+B\left(\frac{1}{1-z_2}\right)^n}{D\left(\frac{1}{1-z_1}\right)^n+E\left(\frac{1}{1-z_2}\right)^n},\
\ \ 0 \leq n \leq n_{0}-1,\label{exp-ben-pure-n0-}\\
S_{ao}(n_0;\lceil n_0 \rceil )&=&R-C\frac{\sum_{i=n_0}^{\infty}\left[A\left(\frac{1}{1-z_1}\right)^i+B\left(\frac{1}{1-z_2}\right)^i\right]}{\sum_{i=n_0}^{\infty}\left[D\left(\frac{1}{1-z_1}\right)^i+E\left(\frac{1}{1-z_2}\right)^i\right]},
\label{exp-ben-pure-n0}
\end{eqnarray}
where
 $A,\ B,\ D,\ E,\ z_1,\ z_2$ are given by \eqref{A}-\eqref{E} and
\eqref{z12}.
\end{cor}

\begin{remark} Applying Corollary \ref{expected-benefit-pure} for $n_0=0$
yields the expected net benefit  $S_{ao}(0;\lceil 0 \rceil )\equiv
S_{ao}(0;\lfloor \infty \rfloor )$ of an arriving customer, if he
decides to join, when the others follow the `always balk' strategy.
\end{remark}

When the customers follow a $(0,\theta)$-mixed reverse-threshold
strategy, with $\theta\in (0,1)$, we can use the same line of argument with Propositions
\ref{expected-benefit-all-enter} and \ref{expected-benefit-mixed},
using the stationary distribution given by
\eqref{stat-rev-prob-main-0}-\eqref{stat-rev-prob-main-0+}. Then we
have the following Proposition \ref{expected-benefit-reverse}.

\begin{prop}\label{expected-benefit-reverse}
Consider the almost observable model of the stochastic clearing
system in alternating environment, where the customers join the
system according to the $(0,\theta)$-mixed reverse-threshold
strategy for some $\theta\in(0,1)$. Then, the expected net benefit
$S_{ao}(n;\lfloor 0,\theta \rfloor )$ of an arriving costumer, if
he decides to join, given that he finds $n$ customers in the system,
is given by
\begin{eqnarray}
S_{ao}(n;\lfloor 0,\theta \rfloor )&=&
R-C\frac{\sum_{i=n}^{\infty}(1-\theta)^{i-n}\left[A\left(\frac{1}{1-z_1}\right)^i+B\left(\frac{1}{1-z_2}\right)^i\right]}{\sum_{i=n}^{\infty}(1-\theta)^{i-n}\left[D\left(\frac{1}{1-z_1}\right)^i+E\left(\frac{1}{1-z_2}\right)^i\right]},
n\geq 0,\label{exp-ben-reverse-1}
\end{eqnarray}
where
 $A,\ B,\ D,\ E,\ z_1,\ z_2$ are given by \eqref{A}-\eqref{E} and
\eqref{z12}.
\end{prop}

To express the various formulas reported in Propositions \ref{expected-benefit-all-enter}--\ref{expected-benefit-reverse} and in Corollary \ref{expected-benefit-pure} for
the expected net benefit function in a compact, unified way, we
introduce the functions
\begin{eqnarray}
F(n,\theta)&=&\sum_{i=n}^{\infty} (1-\theta)^{i-n} \left[ (RD-CA)
\left( \frac{1}{1-z_1} \right)^i +(RE-CB) \left( \frac{1}{1-z_2}
\right)^i \right],\label{FunctionF}\\
G(n,\theta)&=&\sum_{i=n}^{\infty} (1-\theta)^{i-n} \left[ D \left(
\frac{1}{1-z_1} \right)^i +E \left( \frac{1}{1-z_2} \right)^i
\right],\; n\geq 0,\; \theta \in [0,1],\label{FunctionG}\\
H^U(n)&=&\frac{F(n,1)}{G(n,1)},\;\;\;
H^L(n)=\frac{F(n,0)}{G(n,0)},\; n\geq 0.\label{FunctionH}
\end{eqnarray}
Then we have
\begin{eqnarray}
S_{ao}(n;\lceil \infty \rceil )\equiv
S_{ao}(n;\lfloor 0 \rfloor )&=&\frac{F(n,1)}{G(n,1)}=H^U(n),\; n\geq 0,\label{S-ao-FG-1}\\
S_{ao}(n;\lceil n_0,\theta \rceil )&=&\frac{F(n,1)}{G(n,1)}=H^U(n),\; 0\leq n\leq n_0-1,\label{S-ao-FG-2}\\
S_{ao}(n_0;\lceil n_0,\theta \rceil )&=&\frac{F(n_0,\theta)}{G(n_0,\theta)},\label{S-ao-FG-3}\\
S_{ao}(n_0+1;\lceil n_0,\theta \rceil )&=&\frac{F(n_0,0)-F(n_0,\theta)}{G(n_0,0)-G(n_0,\theta)},\label{S-ao-FG-4}\\
S_{ao}(n;\lceil n_0 \rceil )&=&\frac{F(n,1)}{G(n,1)}=H^U(n),\; 0\leq n\leq n_0-1,\label{S-ao-FG-5}\\
S_{ao}(n_0;\lceil n_0 \rceil )&=&\frac{F(n_0,0)}{G(n_0,0)}=H^L(n_0),\label{S-ao-FG-6}\\
S_{ao}(0;\lceil 0 \rceil )\equiv S_{ao}(0;\lfloor \infty \rfloor )&=&\frac{F(0,0)}{G(0,0)}=H^L(0),\label{S-ao-FG-7}\\
S_{ao}(n;\lfloor 0,\theta \rfloor )&=&\frac{F(n,\theta)}{G(n,\theta)},\;
n\geq 0.\label{S-ao-FG-8}
\end{eqnarray}

\section{The almost observable case: Equilibrium strategies}\label{Observable-Section-Equilibrium}

As we have already discussed in the beginning of Section
\ref{Observable-Section-Preliminaries}, it seems plausible that
threshold strategies are adopted by the customers when the `fast
service' environmental state coincides with the `less congested'
environmental state, i.e. when
$(\mu_1-\mu_2)(\rho_1-\rho_2)<0$. On the contrary, reverse-threshold strategies are plausible when
the `fast service' environmental state coincides with the `more
congested' environmental state, i.e. when the opposite inequality
holds. This intuitive finding is associated with the monotonicity of
$H^U(n)$ which plays a key role in the subsequent analysis. More
specifically, we have the following Proposition
\ref{monotonicityHUn}.

\begin{prop}\label{monotonicityHUn}
We have the following equivalences:
\begin{eqnarray}
&&H^U(n) \mbox{ is strictly decreasing} \Leftrightarrow AE-BD>0
\Leftrightarrow
(\mu_1-\mu_2)(\rho_1-\rho_2)<
0.\label{HUn-decreasing-condition}\\
&&H^U(n) \mbox{ is constant} \Leftrightarrow AE-BD=0 \Leftrightarrow
\mu_1=\mu_2 \mbox{ or } \rho_1=\rho_2.\label{HUn-constant-condition}\\
&&H^U(n) \mbox{ is strictly increasing} \Leftrightarrow AE-BD<0
\Leftrightarrow
(\mu_1-\mu_2)(\rho_1-\rho_2)>
0.\label{HUn-increasing-condition}
\end{eqnarray}
\end{prop}

The proof of this proposition is omitted, since its first case
follows easily by simple algebraic manipulations that start from
the relation $H^U(n+1)-H^U(n)<0$ and lead to $AE-BD>0$ and
$(\mu_1-\mu_2)(\rho_1-\rho_2)<
0$, through successive equivalences. The other two cases are
treated similarly. Moreover, the monotonicity of the function $\frac{F(n,\theta)}{G(n,\theta)}$ with
respect to $\theta$ depends on the sign of
$(\mu_1-\mu_2)(\rho_1-\rho_2)$.
Specifically, we have the following Proposition
\ref{monotonicityFG}.

\begin{prop}\label{monotonicityFG}
We have the following equivalences:
\begin{eqnarray}
&&\frac{F(n,\theta)}{G(n,\theta)} \mbox{ is strictly increasing in $\theta$}
\Leftrightarrow AE-BD>0 \Leftrightarrow
(\mu_1-\mu_2)(\rho_1-\rho_2)<
0.\label{FG-increasing-condition}\\
&&\frac{F(n,\theta)}{G(n,\theta)} \mbox{ is constant in $\theta$}
\Leftrightarrow AE-BD=0 \Leftrightarrow
\mu_1=\mu_2 \mbox{ or } \rho_1=\rho_2.\label{FG-constant-condition}\\
&&\frac{F(n,\theta)}{G(n,\theta)} \mbox{ is strictly decreasing in $\theta$}
\Leftrightarrow AE-BD<0 \Leftrightarrow
(\mu_1-\mu_2)(\rho_1-\rho_2)>
0.\label{FG-decreasing-condition}
\end{eqnarray}
\end{prop}
The proof of this proposition is also omitted, since the result
is deduced easily after some algebra. We now state some properties of $F(n,\theta)$, $G(n,\theta)$ and
$H^U(n)$, $H^L(n)$ that we will use in the sequel. Their proof is
straightforward from their definition and thus it is omitted.

\begin{lem} The functions $F(n,\theta)$, $G(n,\theta)$ satisfy the
following properties:
\begin{eqnarray}
&&F(n,\theta)=\sum_{i=n}^{\infty} (1-\theta)^{i-n} F(i,1)=F(n,1)+(1-\theta)F(n+1,\theta),\;
n\geq 0,\;
\theta\in[0,1],\label{PropF2}\\
&&G(n,\theta)=\sum_{i=n}^{\infty} (1-\theta)^{i-n}
G(i,1)=G(n,1)+(1-\theta)G(n+1,\theta),\; n\geq 0,\;
\theta\in[0,1],\label{PropG1}\\
&&G(n,\theta)>0,\; n\geq 0,\; \theta\in[0,1],\label{PropG2}\\
&&G(n,\theta) \mbox{ is strictly decreasing with respect to $\theta$
for any fixed $n\geq 0$}.\label{PropG3}
\end{eqnarray}
\end{lem}

Note that properties \eqref{PropG2} and \eqref{PropG3} of $G(n,\theta)$ assure that all denominators in \eqref{S-ao-FG-1}-\eqref{S-ao-FG-8} are positive.

The intuitive discussion at the beginning of Section \ref{Observable-Section-Preliminaries} in combination with Propositions \ref{monotonicityHUn} and \ref{monotonicityFG} suggests that we should methodologically proceed by considering separately three cases, corresponding to the sign (negative, positive or zero) of $(\mu_1-\mu_2)(\rho_1-\rho_2)$.

\subsection{Case A: $(\mu_1-\mu_2)(\rho_1-\rho_2)<0$}\label{case-a}

In Case A, we will prove that an equilibrium threshold
strategy always exists. Moreover, we will present a systematic procedure for determining all equilibrium threshold strategies. We first introduce several quantities that we will need in the sequel.

\begin{defn} Suppose that
\begin{equation}
(\mu_1-\mu_2)(\rho_1-\rho_2)<0.\label{Condition-threshold-policies}
\end{equation}
We define
\begin{eqnarray}
n_U&=&\inf \{ n\geq 0 : F(n,1)<0\},\label{defn-nu}\\
n_L&=&\inf \{ n\geq 0 : F(n,0)\leq 0\}, \label{defn-nl}\\
n_U^-&=& \inf \{ n \geq 0 : F(n,1)\leq 0\}, \label{defn-nu-}\\
n_L^+&=& \inf \{ n \geq 0 : F(n,0) <0\}. \label{defn-nl+}
\end{eqnarray}
\end{defn}

Then, we have several properties of $n_U$, $n_L$, $n_U^-$ and $n_L^+$ that we summarize in the following Lemma \ref{HUn-decreasing-3cases}.

\begin{lem}\label{HUn-decreasing-3cases}
Suppose that \eqref{Condition-threshold-policies} holds. Then, there are
three cases:
\begin{description}
\item [Case I:] $H^U(0)< 0$.\\
Then
\begin{eqnarray}
&&n_U=n_L=n_U^-=n_L^+=0, \label{all-n-0}\\
&&F(n,\theta)<0,\ \ n\geq 0,\ \ \theta\in [0,1],\label{decreasing-case1-2}\\
&&F(n,0)-F(n,\theta)<0,\ \ n\geq 0,\ \ \theta\in
(0,1].\label{decreasing-case1-3}
\end{eqnarray}
\item [Case II:] $H^U(0)\geq 0$ and
$\lim_{n\rightarrow\infty} H^U(n)< 0$.\\
Then
\begin{eqnarray}
&&1\leq n_U<\infty, \label{nu-proper}\\
&&F(n,1)>0,\ \ 0\leq n\leq n_U-2,\label{decreasing-case2-1}\\
&&F(n_U-1,1)\geq 0,\label{decreasing-case2-2}\\
&&F(n,1)<0,\ \ n\geq n_U.\label{decreasing-case2-3}
\end{eqnarray}
and
\begin{eqnarray}
&&0\leq n_L \leq n_U,\label{nl-proper}\\
&&F(n,0)>0,\ \ 0\leq n\leq n_L-1,\label{decreasing-case2-4}\\
&&F(n_L,0)\leq 0,\label{decreasing-case2-5}\\
&&F(n,0)<0,\ \ n\geq n_L+1.\label{decreasing-case2-6}
\end{eqnarray}
Moreover,\\
\begin{eqnarray}
&& n_L^+ = \left\{ \begin{array}{ll}
n_L, & \textrm{if $F(n_L,0)< 0$}\\
n_L+1, & \textrm{if $F(n_L,0)= 0,$}
\end{array} \right.  \label{decreasing-case2-7}\\
&&n_U^- = \left\{ \begin{array}{ll}
n_U, & \textrm{if $F(n_U-1,1)> 0$}\\
n_U-1, & \textrm{if $F(n_U-1,1)= 0.$}
\end{array} \right. \label{decreasing-case2-8}
\end{eqnarray}
For every $n_0\in
\{n_L^+,\ldots,n_U^- -1 \}$, a unique solution $\theta (n_0)\in (0,1)$ of the equation $F(n_0,\theta)=0$ exists with respect to $\theta$, i.e.
\begin{equation}
F(n_0,\theta(n_0))=0,\; n_L^+ \leq n \leq n_U^- -1.\label{decreasing-case2-9}
\end{equation}

\item [Case III:] $\lim_{n\rightarrow\infty} H^U(n)\geq 0$.\\
Then
\begin{eqnarray}
&&n_U=n_L=n_U^-=n_L^+=\infty, \label{all-n-infty}\\
&&F(n,\theta)>0,\ \ n\geq 0,\ \ \theta\in [0,1],\label{decreasing-case3-2}\\
&&F(n,0)-F(n,\theta)>0,\ \ n\geq 0,\ \ \theta\in
(0,1].\label{decreasing-case3-3}
\end{eqnarray}
\end{description}
\end{lem}

\begin{prooof}{Proof} In Case I, the condition $H^U(0)< 0$ in combination with the monotonicity of $H^U(n)$ (due to \eqref{HUn-decreasing-condition}) implies that $H^U(n)<0$, $n\geq 0$. Then, we have that $F(n,1)<0$, $n\geq 0$ and therefore
$F(n,\theta)=\sum_{i=n}^{\infty}(1-\theta)^{i-n}F(i,1)<0,\ n\geq
0,\ \theta\in [0,1]$ and
$F(n,0)-F(n,\theta)=\sum_{i=n}^{\infty}[1-(1-\theta)^{i-n}]F(i,1)<0,\
n\geq 0,\ \theta\in (0,1]$.

In Case II, the conditions $H^U(0)\geq 0$ and
$\lim_{n\rightarrow\infty} H^U(n)< 0$, in combination with the
condition \eqref{HUn-decreasing-condition} for the monotonicity of $H^U(n)$ imply \eqref{nu-proper}-\eqref{decreasing-case2-3}.

Equation \eqref{decreasing-case2-3} implies that $F(n_U,0)=\sum_{i=n_U}^{\infty} F(i,1)<0$ and we conclude \eqref{nl-proper}.
Note also that by the definition of $n_L$ we have immediately \eqref{decreasing-case2-4}-\eqref{decreasing-case2-5}. Moreover, we have that $F(n,0)=\sum_{i=n}^{\infty}F(i,1)<0$, for $n \geq n_U$. For $n$ with $n_L+1\leq n \leq n_U-1$ we have also $F(n,0)<0$. Indeed, suppose that there exists an $n$ with $n_L+1\leq n \leq n_U-1$ such that $F(n,0)\geq 0$. Then we would have $F(n-1,0)=F(n-1,1)+F(n,0)>0$ and inductively we would obtain that $F(n_L,0)>0$, a contradiction because of \eqref{decreasing-case2-5}. Thus $F(n,0)<0$, for all $n\geq n_L+1$ and we obtain \eqref{decreasing-case2-6}.

Equations \eqref{decreasing-case2-7} and \eqref{decreasing-case2-8} are immediate from \eqref{decreasing-case2-1}-\eqref{decreasing-case2-3} and \eqref{decreasing-case2-4}-\eqref{decreasing-case2-6} respectively. Consider, now, an $n_0\in \{n_L^+,\ldots,n_U^- -1 \}$. Then we have that $\frac{F(n_0,1)}{G(n_0,1)}>0$ (since $n_0\leq n_U^- -1$ - see \eqref{decreasing-case2-1}) and $\frac{F(n_0,0)}{G(n_0,0)}<0$ (since $n_0\geq n_L^+$ - see \eqref{decreasing-case2-6}-\eqref{decreasing-case2-7}).  By condition
\eqref{FG-increasing-condition}, we have that $\frac{F(n_0,\theta)}{G(n_0,\theta)}$
is a strictly increasing and continuous function of $\theta$, so by Bolzano's theorem we conclude that there exists a unique solution $\theta(n_0)\in(0,1)$ of the equation $\frac{F(n_0,\theta)}{G(n_0,\theta)}=0$. Thus, we obtain \eqref{decreasing-case2-9}.

In Case III, the condition $\lim_{n\rightarrow\infty} H^U(n)\geq 0$ in combination with the condition \eqref{HUn-decreasing-condition} for the monotonicity of $H^U(n)$ implies that $H^U(n)>0$, $n\geq 0$ which gives $F(n,1)>0$, $n\geq 0$ and $F(n,\theta)=\sum_{i=n}^{\infty}(1-\theta)^{i-n}F(i,1)>0$, for $n \geq 0$ and $\theta\in [0,1]$.  Moreover,
$F(n,0)-F(n,\theta)=\sum_{i=n}^{\infty}[1-(1-\theta)^{i-n}]F(i,1)>0$, $n\geq 0$, $\theta\in (0,1]$. Thus, we conclude \eqref{all-n-infty}-\eqref{decreasing-case3-3}.
\end{prooof}\\

Using Lemma \ref{HUn-decreasing-3cases} we will now
prove the existence of threshold equilibrium strategies, when \eqref{Condition-threshold-policies} holds. We present the results in the following Theorem \ref{decreasing}.

\begin{thm}\label{decreasing}
In the almost observable model of the stochastic clearing system
in alternating environment, where \eqref{Condition-threshold-policies}
holds, equilibrium threshold strategies always exist. In particular, in the three cases of Lemma \ref{HUn-decreasing-3cases} we have:
\begin{description}
\item [Case I:] $H^U(0)< 0$.\\
 Then, there is a unique equilibrium threshold strategy, the $\lceil 0 \rceil $-strategy (always to balk).
\item [Case II:] $H^U(0)\geq 0$ and
$\lim_{n\rightarrow\infty} H^U(n)< 0$.\\
Then, an equilibrium pure threshold strategy always exists. Moreover, the equilibrium strategies within the class of all pure threshold strategies are the strategies $\lceil n_0 \rceil $ with $n_0=n_L,n_L+1,\ldots,n_U$. Also, the equilibrium strategies within the class of genuinely mixed threshold strategies are the strategies $\lceil n_0,\theta(n_0) \rceil $ with $n_0\in \{n_L^+,\ldots,n_U^- -1 \}$ and $\theta(n_0)$ the unique solution in $(0,1)$ of $F(n_0,\theta)=0$ with respect to $\theta$.
\item [Case III:] $\lim_{n\rightarrow\infty}
H^U(n)\geq 0$.\\\
Then, there is a unique equilibrium threshold strategy, the $\lceil \infty \rceil $-strategy (always to join).
\end{description}
\end{thm}

\begin{prooof}{Proof}
Case I: Consider a tagged customer at his arrival instant and
assume that all other customers follow an
$\lceil n_0 \rceil $ strategy for some $n_0\geq 0$. Inequality
\eqref{decreasing-case1-2} and
relations \eqref{S-ao-FG-5} and \eqref{S-ao-FG-6} imply that the
expected net benefit of the tagged customer, when he
finds $n$ customers and decides to join is
$S_{ao}(n;\lceil n_0 \rceil )<0$, for $0\leq n\leq n_0$. Thus, he always prefers to balk and his best response against $\lceil n_0 \rceil $ is $\lceil 0 \rceil $.

We now assume that all other customers follow an
$\lceil n_0,\theta_0 \rceil $ strategy, for some $n_0\geq 0$ and $\theta_0 \in (0,1)$. Then, if the tagged customer finds $n$ customers at his arrival instant and decides to join, his expected net benefit will be $S_{ao}(n;\lceil n_0,\theta_0 \rceil )<0$ for $0\leq n\leq n_0+1$, from
\eqref{decreasing-case1-2}-\eqref{decreasing-case1-3} and
\eqref{S-ao-FG-2}-\eqref{S-ao-FG-4}. Therefore,
the tagged customer is always unwilling to join and we have that his best response against $\lceil n_0,\theta_0 \rceil $ is $\lceil 0 \rceil $.

If all customers follow the $\lceil \infty \rceil $ strategy,
\eqref{decreasing-case1-2} and \eqref{S-ao-FG-1} yield
$S_{ao}(n;\lceil \infty \rceil )<0$ for $n\geq 0$. Again, due to the negative expected net benefit, it is preferable for the tagged customer to balk. So, his best response against $\lceil \infty \rceil $ is $\lceil 0 \rceil $. Thus, we conclude that the only best response against itself within the class of (pure and mixed) threshold strategies is $\lceil 0 \rceil $.\\

Case II: Consider a tagged arriving customer and suppose that all other customers follow an $\lceil n_0 \rceil $ strategy, for some $n_0\leq n_L-1$. If the tagged customer finds $n_0$ customers and decides to join, his expected net benefit will be $S_{ao}(n_0;\lceil n_0 \rceil )>0$, from \eqref{decreasing-case2-4} and \eqref{S-ao-FG-6}. This implies that when he finds $n_0$ customers, he is willing to join. Thus,
$\lceil n_0 \rceil $ cannot be a best response against itself. So such a strategy cannot be an equilibrium.

Consider, now, a tagged arriving customer and suppose that all other customers follow an $\lceil n_0 \rceil $ strategy, for some $n_0\geq n_U+1$. Using \eqref{S-ao-FG-5} and \eqref{decreasing-case2-3}, we have that
$S_{ao}(n;\lceil n_0 \rceil )<0$, for $n_U\leq n\leq n_0-1$. This means that when the tagged customer finds $n$ customers, with $n_U\leq n\leq n_0-1$, then he is unwilling to enter. Thus, the $\lceil n_0 \rceil $ strategy cannot be an equilibrium. We conclude that the search for equilibrium strategies within the class of pure threshold strategies should be restricted to strategies $\lceil n_0 \rceil $ with $n_L\leq n_0\leq n_U$.

We mark an arriving customer and we assume that all other customers follow an $\lceil n_0 \rceil $ strategy, for some $n_0$ with $n_L\leq n_0\leq n_U$. From \eqref{S-ao-FG-5}, \eqref{S-ao-FG-6}, \eqref{decreasing-case2-1}, \eqref{decreasing-case2-2}, \eqref{decreasing-case2-5} and \eqref{decreasing-case2-6}, we have that the
expected net benefit of a customer who finds $n$ customers upon arrival and decides to join is $S_{ao}(n;\lceil n_0 \rceil )\geq 0$, for $0\leq n\leq n_0-1$ and $S_{ao}(n_0;\lceil n_0 \rceil )\leq 0$. Thus $\lceil n_0 \rceil $ is a best response against itself and we conclude that all such strategies are equilibrium strategies.

To finish with our search for equilibrium strategies in the class of pure threshold strategies, we have to
examine the $\lceil \infty \rceil $ strategy. This cannot be an
equilibrium, since \eqref{S-ao-FG-1} and
\eqref{decreasing-case2-3} imply that
$S_{ao}(n;\lceil \infty \rceil )< 0$, for $n\geq n_U$, which means that it is not optimal for the tagged customer to join when he sees $n$ customers for some $n\geq n_U$. Therefore, we conclude that the equilibrium strategies within the class of pure threshold strategies are exactly the strategies $\lceil n_0 \rceil $ with $n_L\leq n_0\leq n_U$.

We will now search for equilibrium strategies in the class of genuinely mixed threshold strategies, i.e. among strategies $\lceil n_0,\theta_0 \rceil $ with $\theta_0\in(0,1)$. A mixed
threshold strategy $\lceil n_0,\theta_0 \rceil $ is an equilibrium if and only if the relations $F(n,1)\geq 0$, for $0\leq n\leq n_0-1$, $F(n_0,\theta_0)=0$ and $F(n_0,0)-F(n_0,\theta_0)\leq 0$ hold (see \eqref{S-ao-FG-2}-\eqref{S-ao-FG-4}). A moment of reflection shows that $\lceil n_0,\theta_0 \rceil $ may be an equilibrium only if $\lceil n_0 \rceil $ is an equilibrium (see \eqref{S-ao-FG-2}-\eqref{S-ao-FG-4} in comparison with \eqref{S-ao-FG-5}-\eqref{S-ao-FG-6}). Thus, we should restrict our search for equilibrium genuinely mixed threshold strategies to strategies $\lceil n_0,\theta_0 \rceil $ with $n_0=\ n_L,\ n_L+1,\ \ldots,\ n_U$.

If $F(n_L,0)=0$, then there does not exist $\theta \in (0,1)$ such that $F(n_L,\theta)=0$, since $\frac{F(n,\theta)}{G(n,\theta)}$ is strictly decreasing. Therefore, $\lceil n_L,\theta \rceil $ cannot be an
equilibrium strategy for any $\theta \in (0,1)$. Similarly, if $F(n_U-1,1)=0$, then the strategy $\lceil n_U-1,\theta \rceil $ cannot be an
equilibrium strategy for any $\theta \in (0,1)$. Moreover,  $\lceil n_U,\theta \rceil $ cannot be equilibrium for any $\theta \in (0,1)$, since $F(n_U,\theta)<0$, $\theta \in (0,1)$. Therefore, a strategy $\lceil n_0,\theta_0 \rceil $ with $\theta\in(0,1)$ may be an equilibrium only if $n_L^+ \leq n_0 \leq n_U^- -1$.

Now, for every $n_0\in \{n_L^+,\ldots,n_U^- -1\}$, the only $\lceil n_0,\theta_0 \rceil $ strategy that can be an equilibrium is the one that corresponds to $\theta_0=\theta(n_0)$, since
$F(n_0,\theta(n_0))=0$. Indeed, if all customers follow the
$\lceil n_0,\theta(n_0) \rceil $ strategy, the expected net benefit
for a tagged customer, who finds $n$ other customers and decides to join the system, is $S_{ao}(n;\lceil n_0,\theta(n_0) \rceil )>0$, for $0\leq n\leq n_0-1$, $S_{ao}(n_0;\lceil n_0,\theta(n_0) \rceil )=0$ and
$S_{ao}(n_0+1;\lceil n_0,\theta(n_0) \rceil )<0$, from
\eqref{decreasing-case2-1}, \eqref{decreasing-case2-6} and
\eqref{decreasing-case2-9}. Thus, $\lceil n_0,\theta(n_0) \rceil $ is an equilibrium strategy.\\

Case III: Following the same line of argument as in case I, we now find that when all customers follow a pure threshold strategy $\lceil n_0 \rceil $ or a mixed threshold strategy $\lceil n_0,\theta_0 \rceil $ the expected net benefit function is  always positive. Thus, the best response of a customer is always to join the system. Thus, the only best response against itself in the class of threshold strategies is the $\lceil \infty \rceil $ strategy.
\end{prooof}\\

Note that although pure threshold strategies always exist, it is possible that genuinely mixed threshold strategies do not. This happens if $n_U^- -1<n_L^+$.\\

\subsection{Case B: $(\mu_1-\mu_2)(\rho_1-\rho_2)>0$}\label{case-b}

In Case B, we seek for equilibrium strategies in
the class of reverse-threshold strategies. We will exclude  strategies $\lfloor n_0 \rfloor $ and $\lfloor n_0,\theta_{0} \rfloor $ with $n_{0}\geq 1$. Indeed, all these strategies prescribe to balk, when a tagged arriving customer sees an empty system. Thus, under such a strategy, the system remains continuously empty, after the first service completion. Therefore, in steady state, these strategies are equivalent to the `always balk' strategy $\lfloor \infty \rfloor $. Thus, we seek for equilibrium strategies only in the set $\mathcal{S}_{r-t}=\{\lfloor 0 \rfloor ,\
\lfloor \infty \rfloor \}\cup \{\lfloor 0,\theta_{0} \rfloor :
\theta_{0}\in(0,1)\}$. We first introduce several quantities that we will use in the sequel.

\begin{defn}\label{definition-m}
Suppose that
\begin{equation}
(\mu_1-\mu_2)(\rho_1-\rho_2)>
0.\label{Condition-reversethreshold-policies}
\end{equation}
We define
\begin{eqnarray}
m_U&=&\inf \{ n\geq 0 : F(n,1)>0\},\label{defn-mu}\\
m_L&=&\inf \{ n\geq 0 : F(n,0)\geq 0\}, \label{defn-ml}\\
m_U^-&=& \inf \{ n \geq 0 : F(n,1)\geq 0\}, \label{defn-mu-}\\
m_L^+&=& \inf \{ n \geq 0 : F(n,0) >0\}. \label{defn-ml+}
\end{eqnarray}
\end{defn}

Then, we have several properties of $m_U$, $m_L$, $m_U^-$ and $m_L^+$ that we summarize in the following Lemma
\ref{HUn-increasing-3cases}.

\begin{lem}\label{HUn-increasing-3cases}
Suppose that \eqref{Condition-reversethreshold-policies} holds. Then, there are three cases:
\begin{description}
\item [Case I:] $H^U(0)> 0$.

Then
\begin{eqnarray}
&&m_U=m_L=m_U^-=m_L^+=0,\label{all-m-0}\\
&&F(n,\theta)>0,\ \ n\geq 0,\ \ \theta\in
[0,1].\label{increasing-case1-2}
\end{eqnarray}

\item [Case II:]
$H^U(0)\leq 0$ and
$\lim_{n\rightarrow\infty} H^U(n)> 0$.

Then
\begin{eqnarray}
&&1\leq m_U < \infty,\label{mu-proper}\\
&&F(n,1)<0,\ \ 0\leq n\leq m_U-2,\label{increasing-case2-1}\\
&&F(m_U-1,1)\leq 0,\label{increasing-case2-2}\\
&&F(n,1)>0,\ \ n\geq m_U.\label{increasing-case2-3}
\end{eqnarray}
and
\begin{eqnarray}
&&0\leq m_L \leq m_U,\label{ml-proper}\\
&&F(n,0)<0,\ \ 0\leq n\leq m_L-1,\label{increasing-case2-4}\\
&&F(m_L,0)\geq 0,\label{increasing-case2-5}\\
&&F(n,0)>0,\ \ n\geq m_L+1.\label{increasing-case2-6}
\end{eqnarray}
Moreover,
\begin{eqnarray}
&& m_L^+ = \left\{ \begin{array}{ll}
m_L, & \textrm{if $F(m_L,0)> 0$}\\
m_L+1, & \textrm{if $F(m_L,0)= 0$,}
\end{array} \right.\label{increasing-case2-7}\\
&& m_U^- = \left\{ \begin{array}{ll}
m_U, & \textrm{if $F(m_U-1,1)< 0$}\\
m_U-1, & \textrm{if $F(m_U-1,1)= 0$.}
\end{array} \right.\label{increasing-case2-8}
\end{eqnarray}
If $m_L^+=0$ and $m_U^-\geq 1$, then there exists a unique
$\theta(0)\in \ (0,1)$ such that \begin{eqnarray}
&&F(0,\theta(0))=0,\label{increasing-case2-9}\\
&&F(n,\theta(0))> 0,\ n\geq 1.\label{increasing-case2-10}
\end{eqnarray}

\item [Case III:] $\lim_{n\rightarrow\infty} H^U(n)\leq 0$.\\
Then
\begin{eqnarray}
&&m_U=m_L=m_U^-=m_L^+=\infty,\label{all-m-infty}\\
&&F(n,\theta)<0,\ n\geq 0,\; \theta\in [0,1].\label{increasing-case3-2}
\end{eqnarray}
\end{description}
\end{lem}

We omit the proof of Lemma \ref{HUn-increasing-3cases} as it is completely analogous to the proof of Lemma \ref{HUn-decreasing-3cases}. We are now in position to prove the existence and uniqueness of reverse-threshold strategies, when \eqref{Condition-reversethreshold-policies} holds. We present the results in the following Theorem \ref{increasing}. The statements about the uniqueness of the reverse-threshold equilibrium strategies should be interpreted within the class $\mathcal{S}_{r-t}=\{\lfloor 0 \rfloor ,\
\lfloor \infty \rfloor \}\cup \{\lfloor 0,\theta_{0} \rfloor :
\theta_{0}\in(0,1)\}$ of reverse-threshold strategies.

\begin{thm}\label{increasing}
In the almost observable model of the stochastic clearing system
in alternating environment, where \eqref{Condition-reversethreshold-policies}
holds, equilibrium reverse-threshold strategies always exist. In particular, in the three cases of Lemma \ref{HUn-increasing-3cases} we have:
\begin{description}
\item [Case I:] $H^U(0)> 0$.

Then, there is a unique equilibrium reverse-threshold strategy, the $\lfloor 0 \rfloor $ strategy (`always to join').
\item [Case 2.] $H^U(0)\leq 0$ and
$\lim_{n\rightarrow\infty} H^U(n)> 0$.

If $m_U^- =0$, the $\lfloor 0 \rfloor $ strategy (`always to join') is the unique equilibrium reverse-threshold strategy. If $m_L^+\geq 1$, then the $\lfloor \infty \rfloor $ strategy (`always to balk') is the unique equilibrium reverse-threshold strategy. Otherwise, the $\lfloor 0,\theta(0) \rfloor $ strategy is the unique equilibrium reverse-threshold strategy.

\item [Case III:] $\lim_{n\rightarrow\infty} H^U(n)\leq 0$.

Then, there is a unique equilibrium reverse-threshold strategy, the $\lfloor \infty \rfloor $ strategy (`always to balk').
\end{description}
\end{thm}

\begin{prooof}{Proof}
Case I: Consider a tagged customer at his arrival instant and
assume that all other customers follow the
$\lfloor 0 \rfloor $ strategy. Inequality \eqref{increasing-case1-2}
and relation \eqref{S-ao-FG-1} imply that his expected net
benefit, when he finds $n$ customers and decides to join is
$S_{ao}(n;\lfloor 0 \rfloor )>0$, for $n\geq 0$. Thus, he always prefers to join so his best response against $\lfloor 0 \rfloor $ is $\lfloor 0 \rfloor $ itself.

Similarly, let mark an arriving customer and suppose that all other customers follow a $\lfloor 0,\theta_0 \rfloor $ strategy, for some $\theta_0 \in (0,1)$. Then, the expected net benefit of the tagged customer, who finds $n$ customers at his arrival instant and decides to join, will be $S_{ao}(n;\lfloor 0,\theta_0 \rfloor )> 0$, for $n\geq 0$ due to \eqref{increasing-case1-2} and \eqref{S-ao-FG-8}. Therefore, the tagged customer is always willing to join and we have that $\lfloor 0 \rfloor $ is the best response against $\lfloor 0,\theta_0 \rfloor $.

If all customers follow the $\lfloor \infty \rfloor $ strategy, equations \eqref{increasing-case1-2} and \eqref{S-ao-FG-7} imply that $S_{ao}(0;\lfloor \infty \rfloor )> 0$, so the tagged customer prefers to join. Thus, we have again that $\lfloor 0 \rfloor $ is the best response against $\lfloor \infty \rfloor $. So the only reverse-threshold strategy which is best response against itself is the $\lfloor 0 \rfloor $ strategy.

Case II: Assume that $m_U^- =0$. Then $F(0,1)=0$ and
$m_U=1$. Consider now a tagged customer at his arrival instant and suppose that all other customers follow the
$\lfloor 0 \rfloor $ strategy. Inequality \eqref{increasing-case2-3}
and relation \eqref{S-ao-FG-1} imply that his expected net benefit, when he finds $n$ customers and decides to join is
$S_{ao}(n;\lfloor 0 \rfloor )\geq 0$, for $n\geq 0$. Thus, $\lfloor 0 \rfloor $ is a best response to itself.

Assume, now, that $m_L^+\geq 1$, which means that $F(0,0)\leq 0$. If we consider a tagged arriving customer and suppose that
all other customers follow the $\lfloor \infty \rfloor $ strategy, then the tagged customer, if he finds $0$ customers and decides to join, has expected net benefit $S_{ao}(0;\lfloor \infty \rfloor )\leq 0$, due to \eqref{S-ao-FG-7}. Thus the strategy $\lfloor \infty \rfloor $ is best response against itself, i.e. it is equilibrium strategy. Otherwise, we will have  $m_L^+=0$. Consider again a tagged customer and suppose that the other customers follow the $\lfloor 0,\theta(0) \rfloor $ strategy. Then, if the tagged customer finds $n$ customers at his
arrival instant and decides to join, his expected net benefit
will be either $S_{ao}(0;\lfloor 0,\theta(0) \rfloor )=0$ if $n=0$, or $S_{ao}(n;\lfloor 0,\theta(0) \rfloor )> 0$, if $n\geq 1$, due to \eqref{increasing-case2-9}-\eqref{increasing-case2-10} and
\eqref{S-ao-FG-8}. Therefore, the $\lfloor 0,\theta(0) \rfloor $ strategy is equilibrium strategy.\\

Case III: Following the same line of argument as in
case I, we now conclude that the expected net benefit function is negative. Thus, the best response to every reverse-threshold strategy is $\lfloor \infty \rfloor $. Thus the only equilibrium reverse-threshold strategy is $\lfloor \infty \rfloor $.
\end{prooof}\\

\subsection{Case C: $(\mu_1-\mu_2)(\rho_1-\rho_2)=0$}\label{case-c}

Case C occurs when $\mu_1=\mu_2$ or $\frac{\lambda_1}{\mu_1}=\frac{\lambda_2}{\mu_2}$. In this case, the distinction `fast environmental state' and `slow environmental state' has no sense or the distinction `more congested environmental state' and `less congested environmental state' has no sense. Therefore, we conclude that the information on the number of customers in the system, does not affect the decision of a tagged arriving customer. A similar analysis is possible as in the other two cases and we have the following Theorem \ref{constant}.

\begin{thm}\label{constant}
In the almost observable model of the stochastic clearing system in alternating environment, where
\begin{equation}
\mu_1=\mu_2 \mbox{ or } \rho_1=\rho_2, \label{HUn-constant-condition1}
\end{equation}
an equilibrium strategy exists within the class of threshold and reverse-threshold strategies. In particular we have the following three cases:
\begin{description}
\item [Case I:] $H^U(0)< 0$.

Then, the unique equilibrium strategy in the class of threshold and reverse-threshold strategies is the $\lceil 0 \rceil \equiv \lfloor \infty \rfloor $ strategy (`always to balk').

\item [Case II:] $H^U(0)= 0$.

Then, every strategy in the class of threshold and reverse-threshold strategies is equilibrium strategy.

\item [Case III:] $H^U(0)> 0$.

Then, the unique equilibrium strategy in the class of threshold and reverse-threshold strategies is the $\lceil \infty \rceil \equiv \lfloor 0 \rfloor $ strategy (`always to join').
\end{description}
\end{thm}

\section{Summary and conclusions}\label{Summary-Conclusions}

\indent In this paper we considered the problem of analyzing
customer strategic behavior, in a clearing system in
alternating environment, where customers decide whether to join
the system or balk upon arrival. We identified four cases with respect to
the level of information provided to arriving customers and
derived the equilibrium strategies for each case. It is important
to notice that in each case we identified all equilibrium
strategies within the appropriate class of strategies. Moreover,
in the almost observable case, which is the most interesting one,
Theorems \ref{decreasing}, \ref{increasing} and \ref{constant}
suggest that the equilibrium strategies in the class of threshold
and reverse-threshold strategies are completely characterized by
the signs of the quantities
$(\mu_1-\mu_2)\left(\rho_1-\rho_2\right)$,
$H^U(n)$, $\lim_{n\rightarrow\infty}H^U(n)$ and $H^L(n)$. Thus, we
can easily combine these theorems and develop an algorithm for determining the equilibrium
strategies. We present the algorithm in pseudo-code form in Figure 1. Figure 2
shows schematically the various cases I,II,III when $(\mu_1-\mu_2)\left(\rho_1-\rho_2\right)<0$.\\

\indent We have also to notice that the results in the almost observable case are qualitatively different for the two cases A and B, where $(\mu_1-\mu_2)(\rho_1-\rho_2)$ is negative and positive respectively. Indeed, in case A, there is, in general an interval of thresholds that constitute equilibrium threshold strategies. On the contrary, in case B, there is a unique equilibrium reverse-threshold strategy. These observations correspond to the regimes of Follow-The-Crowd (FTC) and Avoid-The-Crowd (ATC) as defined in Hassin and Haviv (1997, 2003). Indeed, in case A, where $(\mu_1-\mu_2)(\rho_1-\rho_2)<0$ we have that the `fast service' environmental state coincides with the `less congested' environmental state. Then, we can argue as follows, if we want to compare two threshold strategies with thresholds $n$ and $n+1$: If the customers follow a threshold strategy with threshold $n$ and an arriving customer observes $n$ customers in the system, then he deduces that at least $n$ customers arrived since the last clearing epoch. If the customers follow a threshold strategy with threshold $n+1$ and the arriving customer observes $n$ customers, then he deduces that exactly $n$ customers arrived since the last clearing epoch. Thus, in the latter case, the arriving customer has the sense that the system is less congested and therefore the environmental state is most probably the `fast service' one. We conclude that the arriving customer is more willing to enter the system. Therefore, if the customers adopt a higher threshold, an arriving customer tends to follow them in adopting a higher threshold and we have an FTC situation.\\

\indent On the other hand, in case B, where $(\mu_1-\mu_2)(\rho_1-\rho_2)>0$, we have that the `low service' environmental state coincides with the `less congested' environmental state. The usual definition of the ATC situation is not applicable here, since we consider reverse-threshold instead of threshold strategies. Moreover, under any reverse-threshold strategy $\lfloor n, \theta \rfloor$ with $n\geq 1$, the system remains continuously empty after the first visit of the transportation facility and so we have excluded these strategies in our seek for equilibrium strategies. Thus, in case B, we will limit our intuitive discussion of the ATC phenomenon to the class of strategies $\{ \lfloor 0, \theta \rfloor : \theta\in[0,1] \}$, as we have already done in the analysis of subsection \ref{case-b}. Suppose that the customers follow a reverse-threshold strategy $\lfloor n, \theta \rfloor$ and then they move to another reverse-threshold strategy $\lfloor n, \theta' \rfloor$ with $\theta'>\theta$. Consider now an arriving customer that finds $0$ customers in the system.  Knowing the strategies of the other customers, the arriving customer has the sense that the system is in the less congested environmental state in the second case, where the customers enter with probability $\theta'$. Indeed in this case, the customers are more willing to join than in the first case (since $\theta'>\theta$) so the information of an empty system imply that it is more probable that the system is in the less congested environmental state. Therefore, the customer becomes less willing to enter, as the less congested environmental state coincides with the low service state. Thus, when the other customers increase the probability of entering, the tagged customer tends to decrease his probability of entering, i.e. we have an ATC situation.\\

\indent The focus of this work was on equilibrium analysis. On the
other hand, one can think of a situation where a central planner
employs acceptance policies that maximize the social benefit,
under the various levels of information on the system state. It is
easy to see that in the fully unobservable, the fully observable and
the almost unobservable cases the strategies that maximize the social
benefit are the equilibrium strategies. This coincidence between equilibrium and socially optimal strategies
can be explained by the total removals. Since the server removes
all customers at service completion epochs, each customer who
decides to join does not impose any externalities to other
customers. In the almost
observable case equilibrium and socially optimal strategies are
identical except from the case where $H^U(n)$ is strictly
decreasing, $H^U(0)\geq 0$ and $\lim_{n\rightarrow\infty} H^U(n)<
0$. In this case the unique socially optimal strategy is the
$\lceil n_U \rceil $ strategy, which is also an equilibrium.

\section{Bibliography}\label{Bibliography-Section}

\begin{enumerate}
\item Artalejo, J.R. and Gomez-Corral, A. (1998) Analysis
of a stochastic clearing system with repeated attempts.
\textit{Communications in Statistics - Stochastic Models}
\textbf{14}, 623-645.

\item Burnetas, A. and Economou, A. (2007) Equilibrium customer
strategies in a single server Markovian queue with setup times.
\textit{Queueing Systems} \textbf{56}, 213-228.

\item Economou, A. (2003) On the control of a compound immigration
process through total catastrophes. \textit{European Journal of
Operational Research} \textbf{147}, 522-529.

\item Economou, A. and Fakinos, D. (2003) A continuous-time
Markov chain under the influence of a regulating point process and
applications in stochastic models with catastrophes.
\textit{European Journal of Operational Research} \textbf{149},
625-640.

\item Economou, A. and Fakinos, D. (2008) Alternative
approaches for the transient analysis of Markov chains with
catastrophes. \textit{Journal of Statistical Theory and Practice}
\textbf{2}, 183-197.

\item Economou, A. and Kanta, S. (2008a) Optimal balking strategies
and pricing for the single server Markovian queue with compartmented
waiting space. \textit{Queueing Systems} \textbf{59}, 237-269.

\item Economou, A. and Kanta, S. (2008b) Equilibrium balking strategies
in the observable single-server queue with breakdowns and repairs.
\textit{Operations Research Letters} \textbf{36}, 696-699.

\item Edelson, N.M. and Hildebrand, K. (1975) Congestion tolls
for Poisson queueing processes. \textit{Econometrica} \textbf{43},
81-92.

\item Gani, J. and Swift, R.J. (2007) Death and birth-death and
immigration processes with catastrophes. \textit{Journal of
Statistical Theory and Practice} \textbf{1}, 39-48.

\item Guo, P. and Zipkin, P. (2007) Analysis and comparison of queues
with different levels of delay information. \textit{Management
Science} \textbf{53}, 962-970.

\item Hassin, R. (2007) Information and uncertainty in a queuing system.
\textit{Probability in the Engineering and Informational Sciences}
\textbf{21}, 361-380.

\item Hassin, R. and Haviv, M. (1997) Equilibrium threshold strategies:
the case of queues with priorities. \textit{Operations Research}
\textbf{45}, 966-973.

\item Hassin, R. and Haviv, M. (2003) \textit{To Queue or Not to
Queue: Equilibrium Behavior in Queueing Systems.} Kluwer Academic
Publishers, Boston.

\item Kim, K. and Seila, A.F. (1993) A generalized cost model for
stochastic clearing systems. \textit{Computers and Operations
Research} \textbf{20}, 67-82.

\item Kyriakidis, E.G. (1994) Stationary probabilities for a simple
immigration-birth-death process under the influence of total
catastrophes. \textit{Statistics and Probability Letters}
\textbf{20}, 239-240.

\item Kyriakidis, E.G. (1999a) Optimal control of a truncated
general immigration process through total catastrophes.
\textit{Journal of Applied Probability} \textbf{36}, 461-472.

\item Kyriakidis, E.G. (1999b) Characterization of the optimal
policy for the control of a simple immigration process through total
catastrophes. \textit{Operations  Research Letters} \textbf{24},
245-248.

\item Kyriakidis, E.G. and Dimitrakos, T.D. (2005) Computation of
the optimal policy for the control of a compound immigration process
through total catastrophes. \textit{Methodology and Computing in
Applied Probability} \textbf{7}, 97-118.

\item Naor, P. (1969) The regulation of queue size by levying tolls.
\textit{Econometrica} \textbf{37}, 15-24.

\item Serfozo, R. and Stidham, S. (1978) Semi-stationary
clearing processes. \textit{Stochastic Processes and their
Applications} \textbf{6}, 165-178.

\item Stidham, S.Jr. (1974) Stochastic clearing systems.
\textit{Stochastic Processes and their Applications} \textbf{2},
85-113.

\item Stidham, S.Jr. (1977) Cost models for
stochastic clearing systems. \textit{Operations Research}
\textbf{25}, 100-127.

\item Stidham, S.Jr. (2009) \textit{Optimal Design of Queueing
Systems.} CRC Press, Taylor and Francis Group, Boca Raton.

\item Stirzaker, D. (2006) Processes with catastrophes.
\textit{Mathematical Scientist} \textbf{31}, 107-118.

\item Stirzaker, D. (2007) Processes with random regulation.
\textit{Probability in the Engineering and Informational Sciences}
\textbf{21}, 1-17.

\item Sun, W., Guo, P. and Tian, N. (2010) Equilibrium threshold
strategies in observable queueing systems with setup/closedown
times. \textit{Central European Journal of Operations Research}
\textbf{18}, 241-268.

\item Yang, W.S., Kim, J.D. and Chae, K.C. (2002) Analysis of M/G/1
stochastic clearing systems. \textit{Stochastic Analysis and
Applications} \textbf{20}, 1083-1100.

\item Zhang, F. and Wang, J. (2010) Equilibrium analysis of the
observable queue with balking and delayed repairs. \textit{3rd
International Joint Conference on Computational Sciences and
Optimization, CSO 2010: Theoretical Development and Engineering
Practice 2}, art. no. 5533079, 125-129.
\end{enumerate}

\newpage

\textbf{ALGORITHM}\\
\footnotesize{ \textsc{ if $(\mu_1-\mu_2)(\rho_{1}-\rho_{2})<
0$ then\\
\indent if $H^U(0)< 0$ then ``Equilibrium threshold strategies:
$\lceil 0 \rceil$.''\\
\indent elseif $\lim_{n\rightarrow\infty} H^U(n)\geq 0$ then
``Equilibrium threshold strategies: $\lceil \infty \rceil$.''\\
\indent else\\
\indent \indent \% \emph{Compute $n_U:\ n_U=\inf\{n\geq 0:
F(n,1)<0\}$}\\
\indent \indent $n_U=0$\\
\indent \indent while $F(n_U,1)\geq 0$ do\\
\indent \indent \indent $n_U=n_U+1$ \\
\indent \indent endwhile\\
\indent \indent \% \emph{Compute $n_L:\ n_L=\inf\{n\geq 0:
F(n,0)\leq 0\}$}\\
\indent \indent $n_L=n_U$\\
\indent \indent while $F(n_L-1,0)\leq 0$ do\\
\indent \indent \indent $n_L=n_L-1$ \\
\indent \indent endwhile\\
\indent \indent \% \emph{Compute $n_U^-:\ n_U^-=\inf\{n\geq 0:
F(n,1)\leq 0\}$}\\
\indent \indent if $F(n_U-1,1)>0$ then\\
\indent \indent \indent $n_U^-=n_U$\\
\indent \indent else\\
\indent \indent \indent $n_U^-=n_U-1$\\
\indent \indent endif\\
\indent \indent \% \emph{Compute $n_L^+:\ n_L^+=\inf\{n\geq 0:
F(n,0)< 0\}$}\\
\indent \indent if $F(n_L,0)<0$ then\\
\indent \indent \indent $n_L^+=n_L$\\
\indent \indent else\\
\indent \indent \indent $n_L^+=n_L+1$\\
\indent \indent endif\\
\indent \indent ``Equilibrium threshold strategies: $\lceil n_0
\rceil,\
n_0=\ n_L,\ n_L+1,\ ...,\ n_U$.''\\
\indent \indent if $n_L^+\leq n_U^- -1$ then\\
\indent \indent \indent for $n_0=n_L^+:n_U^- -1$\\
\indent \indent \indent \indent Compute
$\theta(n_0):F(n_0,\theta(n_0))=0$\\
\indent \indent \indent endfor\\
\indent \indent \indent ``Equilibrium mixed threshold strategies:
$\lceil n_0,\theta(n_0)\rceil,\ n_0\in \{n_L^+,...,n_U^- -1\}$.''\\
\indent \indent endif\\
\indent endif\\
elseif $(\mu_1-\mu_2)(\rho_{1}-\rho_{2})>
0$ then\\
\indent if $H^U(0)> 0$ then ``Equilibrium reverse-threshold
strategies:
$\lfloor 0 \rfloor$.''\\
\indent elseif $\lim_{n\rightarrow\infty} H^U(n)\leq 0$ then
``Equilibrium reverse-threshold strategies:
$\lfloor \infty \rfloor$.''\\
\indent else\\
\indent \indent if $F(0,1)=0$ then ``Equilibrium reverse-threshold
strategies: $\lfloor 0 \rfloor$.''\\
\indent \indent elseif $F(0,0)\leq 0$ then ``Equilibrium
reverse-threshold
strategies: $\lfloor \infty \rfloor$.''\\
\indent \indent else\\
\indent \indent \indent Compute $\theta(0):F(0,\theta(0))=0$\\
\indent \indent \indent ``Equilibrium reverse-threshold
strategies:
$\lfloor 0,\theta(0) \rfloor$.''\\
\indent \indent endif\\
\indent endif\\
else\\
\indent if $H^U(0)> 0$ then ``Equilibrium threshold strategies:
$\lceil 0 \rceil$.''\\
\indent elseif $H^U(0)< 0$ then ``Equilibrium threshold
strategies:
$\lceil \infty \rceil$.''\\
\indent else ``Equilibrium threshold strategies:$\lceil n \rceil,\ n\geq 0$.''\\
\indent endif\\
endif }}

\normalsize{}
\begin{center}Figure 1: Computation of equilibrium threshold/reverse-threshold strategies \end{center}

 \vspace{0.5cm}
\setlength{\unitlength}{0.6cm}
\begin{picture}(10,7)
\put(0,3){\vector(1,0){8}} \put(1,0){\vector(0,1){6.5}}
\put(7,2.5){$n$}
\put(0.8,2.2){$\diamond$}\put(1.8,1.5){$\diamond$}\put(2.8,1.1){$\diamond$}\put(3.8,0.8){$\diamond$}\put(4.8,0.7){$\diamond$}\put(5.8,0.65){$\diamond$}
\put(0.8,1.4){$\ast$}\put(1.8,1.1){$\ast$}\put(2.8,0.9){$\ast$}\put(3.8,0.7){$\ast$}\put(4.8,0.63){$\ast$}\put(5.8,0.61){$\ast$}
\put(5,5){$\ast:\ H^L(n)$}\put(5,6){$\diamond:\ H^U(n)$}
\put(3,7){Case I} \put(1,-1){Unique equilibrium:}\put(2.5,-2){balk ($\lceil 0 \rceil$)}
\end{picture}
\begin{picture}(10,7)
\put(0,3){\vector(1,0){8}} \put(1,0){\vector(0,1){6.5}}
\put(7,2.5){$n$}
\put(0.8,5.5){$\diamond$}\put(1.8,4.1){$\diamond$}\put(2.8,3.2){$\diamond$}\put(3.8,2.6){$\diamond$}\put(4.8,2.2){$\diamond$}\put(5.8,2){$\diamond$}
\put(0.8,3.5){$\ast$}\put(1.8,2.7){$\ast$}\put(2.8,2.4){$\ast$}\put(3.8,2.2){$\ast$}\put(4.8,2.02){$\ast$}\put(5.8,1.95){$\ast$}
\put(5,5){$\ast:\ H^L(n)$}\put(5,6){$\diamond:\ H^U(n)$}
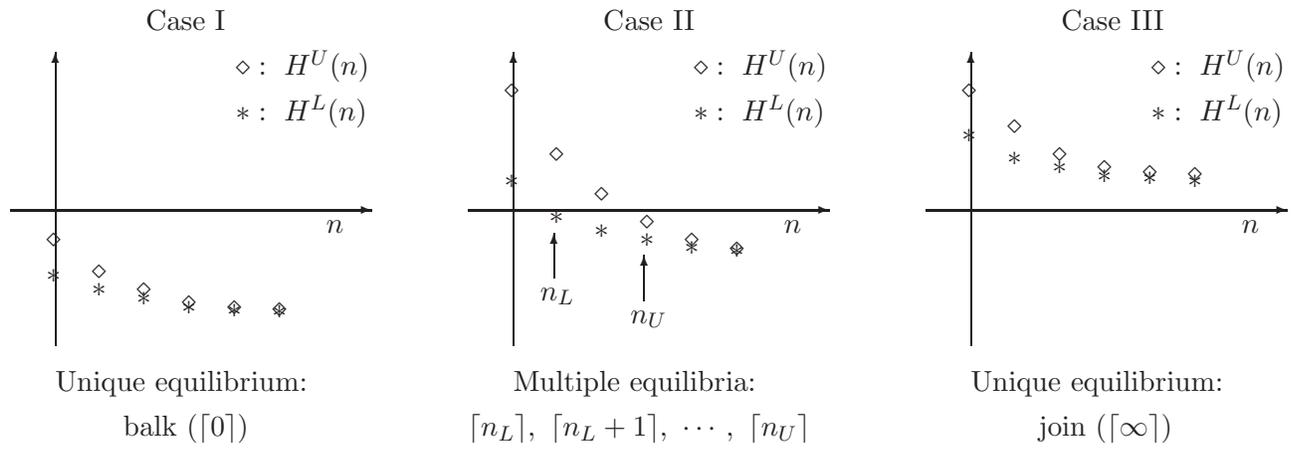
\put(3,7){Case II}\put(1,-1){Multiple equilibria:}\put(0,-2){$\lceil n_L \rceil,\ \lceil n_{L}+1\rceil,\ \cdots,\ \lceil n_U\rceil$}\put(3.9,1){\vector(0,1){1}}\put(3.6,0.5){$n_U$}\put(1.9,1.5){\vector(0,1){1}}\put(1.6,1){$n_L$} \put(-6.5,-4){Figure 2: Case A - $(\mu_1-\mu_2)(\rho_{1}-\rho_{2})<
0$ - Equilibrium threshold strategies }
\end{picture}
\begin{picture}(6,6)
\put(0,3){\vector(1,0){8}} \put(1,0){\vector(0,1){6.5}}
\put(7,2.5){$n$}
\put(0.8,5.5){$\diamond$}\put(1.8,4.7){$\diamond$}\put(2.8,4.1){$\diamond$}\put(3.8,3.8){$\diamond$}\put(4.8,3.7){$\diamond$}\put(5.8,3.65){$\diamond$}
\put(0.8,4.5){$\ast$}\put(1.8,4.0){$\ast$}\put(2.8,3.8){$\ast$}\put(3.8,3.6){$\ast$}\put(4.8,3.55){$\ast$}\put(5.8,3.5){$\ast$}
\put(5,5){$\ast:\ H^L(n)$}\put(5,6){$\diamond:\ H^U(n)$}
\put(3,7){Case III}\put(1,-1){Unique equilibrium:}\put(2.5,-2){join ($\lceil \infty \rceil$)}
\end{picture}\\

\end{document}